\documentclass[letterpaper]{siamltex}
\usepackage[ansinew,latin1]{inputenc}
\usepackage[francais,english]{babel}
\usepackage{amsfonts,amssymb}        %
\usepackage{amsmath}
\usepackage{latexsym}
\usepackage{verbatim}    %
\usepackage{float}
\usepackage{graphicx}
\usepackage{wrapfig}
\usepackage{sidecap}
\usepackage{subfigure}
\def\PID#1{{(#1)}}
\def\BOND{{\mathrm{bond}}}
\def\ANGLE{{\mathrm{angle}}}
\def\TORSION{{\mathrm{torsion}}}
\def\DYN{{\mathrm{dyn}}}
\def\SAMPL{{\mathrm{sampl}}}
\def\RATE{\rho}
\def\ERR{\alpha}

\def\ZETAM{L}
\def\KAPPAM{K}
\def\ID{{\mathrm{Id}\,}}
\def\INTP{\mathrm{int}}
\def\EXTP{\mathrm{e	xt}}

\def\TAG#1{{\qquad ({#1})}}
\def\TAGG#1{{\qquad ({#1})}}

\def\BIGO{\mathcal{O}}
\def\VIZ#1{(\ref{#1})}

\def\WEAKLY{{\rightharpoonup}}

\def\HESS{\mathrm{Hess}\,}
\def\SEP{\,|\,}

\def\DET{\mathrm{det}\,}
\def\TR{\mathrm{Tr}\,}

\def\DIAG{\mathrm{diag}\,}

\def\diss{\M{d}}
\def\nc{n}
\def\SMAN{{\mathcal{M}_{0}}}
\def\SMANZ#1{\mathcal{M}_{#1}}
\def\MASST{M}
\def\MASSTZ{M_z}

\def\nueps{\nu_{\epsilon}}
\def\pen#1{#1_{\nu}}
\def\nusc{\bar{\nu}}
\def\peneps#1{#1_{\nueps}}

\def\DIAG{\mathrm{diag}\,}

\def\Dp{\nabla_{\! p}}
\def\Dq{ \nabla_{\! q} }

\def\Done{ \nabla_{\! 1} }
\def\Dtwo{\nabla_{\! 2}}

\def\penN#1{#1_{\nu_{N}}}

\def\dt{\delta t}
\def\dg{\nabla^{d}}
\def\dl{\Delta_{d}}

\def\EXPECT#1{\E\left[{#1}\right]}

\newcommand{\poisson}[1]{\left\{#1\right\}}

\def\E{\mathbb{E}}

\def\R{\mathbb{R}}

\def\EXP#1{e^{#1}}
\def\ds{\displaystyle}
\def\COMMA{\,,}
\def\PERIOD{\,.}
\def\M#1{{\rm #1 }}

\newcommand{\ph}{\varphi}

\newcommand{\norm}[1]{\left\Vert#1 \right\Vert}

\newcommand{\normop}[1]{|\!|\!|  #1 |\!|\!| }

\newcommand{\abs}[1]{\left\vert#1\right\vert}
\newcommand{\set}[1]{\left\{#1\right\}}

\newcommand{\limop}[1]{\mathop{#1}\limits}

\newenvironment{system} { \left \{
  \begin{alignedat}{2} }
{  \end{alignedat}
  \right.}

\newcommand{\syst}[1]{ \begin{system} #1 \end{system}}

\newcommand{\pare}[1]{\left(#1\right)}

\newenvironment{li} { \left .
  \begin{array}{lcl} }
{  \end{array}
  \right .}
\newcommand{\arr}[1]{ \begin{li} #1 \end{li}}

\newenvironment{tb}{\begin{tabular}{lcl}}{\end{tabular}}

  \newcommand{\bmat}{\begin{pmatrix}}
\newcommand{\emat}{\end{pmatrix}}

\newtheorem{The}{Theorem}[section]
\newtheorem{Lem}[The]{Lemma}
\newtheorem{Pro}[The]{Proposition}

\newtheorem{Def}[The]{Definition}

\newtheorem{scheme}{Scheme}[section]

\newtheorem{Rem}[The]{Remark}
\newtheorem{Exa}[The]{Example}

\newcommand{\bi}{\begin{itemize}}
\newcommand{\ei}{\end{itemize}}

\newcommand{\be}{\begin{equation}}
\newcommand{\ee}{\end{equation}}
\newcommand{\bea}{\begin{eqnarray}}
\newcommand{\eea}{\end{eqnarray}}

\newcommand{\bpro}{\begin{Pro}}
\newcommand{\epro}{\end{Pro}}
\newcommand{\blem}{\begin{Lem}}
\newcommand{\elem}{\end{Lem}}
\newcommand{\bthe}{\begin{The}}
\newcommand{\ethe}{\end{The}}

\newcommand{\bdfn}{\begin{Def}}
\newcommand{\edfn}{\end{Def}}
\newcommand{\brem}{\begin{Rem}}
\newcommand{\erem}{\end{Rem}}
\newcommand{\bexa}{\begin{Exa}}
\newcommand{\eexa}{\end{Exa}}

\title{Implicit Mass-Matrix Penalization of Hamiltonian dynamics with application to exact sampling of stiff systems}
\author{Petr Plech\'a\v{c}\thanks{%
               Department of Mathematics, University of Tennessee, Knoxville, TN 37996-1300 and 
               Oak Ridge National Laboratory, Oak Ridge, TN 37831, USA,
              ({\tt plechac@math.utk.edu}).
	}
	\and
	Mathias Rousset\thanks{%
               INRIA Lille Nord-Europe, Villeneuve d'Ascq, France, ({\tt mathias.rousset@inria.fr})}
}

\begin{document}
\maketitle

\begin{abstract}
An implicit mass-matrix penalization (IMMP) of Hamiltonian dynamics is proposed, and associated dynamical integrators, 
as well as sampling Monte-Carlo schemes, are analyzed for systems with multiple time scales. 
The penalization is based on an extended Hamiltonian with artificial constraints associated with some selected DOFs. 
The penalty parameters enable arbitrary 
tuning of timescales for the selected DOFs.
The IMMP dynamics is shown to be an interpolation between the exact Hamiltonian dynamics 
and the dynamics with rigid constraints. This property translates in the associated 
numerical integrator into a tunable trade-off between stability and dynamical modification. 
Moreover, a penalty that vanishes with the time-step yields order two convergent 
schemes for the exact dynamics. 
Moreover, by construction, the resulting dynamics preserves the canonical equilibrium distribution in position variables, 
up to a computable geometric correcting potential, leading to Metropolis-like unbiased sampling algorithms.
The algorithms can be implemented with a simple modification of standard 
geometric integrators with algebraic constraints imposed on the selected DOFs,
and has no additional complexity in terms of enforcing the constraints and force evaluations.
The properties of the IMMP method are demonstrated numerically on the $N$-alkane model,
showing that the time-step stability region of integrators and the sampling efficiency can be increased 
with a gain that grows with the size of the system. 
This feature is mathematically analyzed for a harmonic atomic chain model.
When a large stiffness parameter is introduced, the IMMP method is shown to 
be asymptotically stable and to
converge towards the heuristically expected Markovian effective dynamics on the slow manifold.
\end{abstract}



\pagestyle{myheadings}
\thispagestyle{plain}
\markboth{P. Plech\'a\v{c}, M. Rousset}%
         {Sampling highly oscillatory systems.}

\section{Introduction}\label{s:intro}
This paper deals with numerical integration and sampling of Hamiltonian
systems with multiple timescales. The main motivation is to develop numerical integration methods for 
the dynamics which resolves certain selected fast degrees of freedom only ''statistically'',
and that can be also used to sample accurately the canonical equilibrium distribution. 
Furthermore, as an ultimate goal, one also seeks good approximation of dynamical 
behavior, at least at large temporal scales.

Hamiltonian systems with multiple timescales typically appear in molecular dynamics (MD)
simulations, which have become, with the aid of increasing computational power, 
a standard tool in many fields of physics, chemistry and biology.
However, extending the simulations to physically relevant time-scales remains
a major challenge for various large molecular systems.
Due to the complexity of implicit methods, the time scales reachable by standard 
numerical methods are usually limited by the rapid oscillations of some particular 
degrees of freedom. Since the sampling dynamics has
to be integrated for long times, the time-step restriction associated with fast oscillations/short
time scales in molecular systems contributes to the high computational cost of such methods.
However, the physical necessity of resolving the fast degrees of freedom in simulations 
is often ambiguous, and efficient
treatment of the fast time scales has motivated new interest in developing numerical schemes
for the integration of such stiff systems.

The problem of integrating stiff forces is relevant both for the direct numerical simulation of the Hamiltonian dynamics, 
as well as for the less restrictive problem of designing a sampling scheme with respect to the canonical ensemble. 
Sampling from the canonical distribution 
can be achieved by Markov chain Monte Carlo (MCMC) algorithms based on a priori
knowledge of possible moves combined with a Metropolis-Hastings acceptance/rejection corrector
(a historical reference is \cite{Met53}). For complex molecular systems, however,
such global moves remain unknown in general, and sampling methods consists generically in using either a Hamiltonian dynamics 
integrator with a thermostat (e.g., a Langevin process),
or its overdamped limit, a drifted random walk (Brownian dynamics)  (see \cite{CanLeg05}
for a review and references on classical sampling methods). Brownian dynamics of systems with multiple time scales suffers 
from similar stability restrictions (see \cite{Fix86,GunBer81}
for some practical issues related to Brownian dynamics simulations in MD).

Broadly speaking, one may start by recognizing two approaches to the numerical treatment of stiff systems:
\begin{description}
\item[{\rm (i)}]   Semi-implicit, multi-step integrators and their variants (e.g., the 
     textbooks
     \cite{LeiRei05, HaiLub02}, or the review paper \cite{JahLub06} and references therein, 
     \cite{Fix86} for Brownian dynamics),
     which attempt to resolve microscopic highly oscillatory dynamical behavior.
\item[{\rm (ii)}] Methods with direct constraints, where the highly oscillatory
     degrees of freedom are constrained to their equilibrium value
     (e.g., \cite{LeiRei05, LeiSke94, GunBer77, Fix78} and references therein).
\end{description}
In spite of their differences the common key feature of all these methods is to balance a trade-off between stability
restrictions and implicit time-stepping form, or in other words,
between the computational effort associated with small time steps, and the computational cost of solving implicit
equations implied by the stiffness.

Although constrained dynamics remove,
in principle, the stiffness of the associated numerical scheme, it introduces new difficulties and numerical problems.
As an approximation to the original dynamics it modifies important features of the system; most importantly,
the original statistical distribution.
The principal goal of the proposed method is to replace direct constraints by implicit mass-matrix
penalization (IMMP), detailed in Section~\ref{s:IMMP},  which integrates fDOFs, but with a tunable mass penalty.
The method designed in this way achieves the two goals:
\begin{description}
 \item[{\rm (i)}]  from the dynamical point of view,  the IMMP method amounts to an 
      appropriate interpolation between exact dynamics and
      constrained dynamics considered in the second family of the methods mentioned above. 
      Moreover, a freely tunable trade-off between dynamical modification and 
      stability is obtained.
\item[{\rm (ii)}] from the sampling point of view, the IMMP dynamics preserves the 
      canonical equilibrium distribution, up to a time step error and an easily 
      computable geometric correcting potential. This leads to 
      Metropolis Monte Carlo methods that sample {\it exactly} the canonical distribution. 
      When using Metropolis schemes, the forces arising from the geometric correcting 
      potential need not be computed.
\end{description}

The idea of adjusting mass tensors in order to slow down fast degrees of freedom goes back to \cite{Ben75}.
In this paper, the author proposes to modify the mass tensor with respect to the Hessian of the potential energy function
in order to confine the frequency spectrum to low frequencies only. Two natural drawbacks of this procedure arise
from the costly computation of the second-order derivatives of the potential, and from the bias introduced when the adjusted mass-tensor
is adapted during the dynamics. Such an approach seems inevitable when the fDOFs are unknown, but in many cases,
the fast degrees of freedom are explicitly given by the structure of the system (e.g., co-valent and angle bonds in molecular chains).
To our knowledge, mass tensor modification have been used in practical MD simulations by increasing the mass of some well-chosen
(e.g., light) atoms \cite{FeeHesBer90, MaoFrie90}. The aim of this paper is to propose a more systematic mass-tensor modification strategy.

The proposed method relies on the assumption that the system Hamiltonian is separable with quadratic kinetic energy
\be\label{e:H}
H(p,q) = \frac{1}{2} p^{T} \MASST^{-1} p  + V(q) \COMMA
\ee
and that the ``fast'' degrees of freedom $(\xi_1,..,\xi_n)$ are explicitly defined,
smooth functions of the system position
\begin{equation}\label{e:xi}
q=(q_1,\dots,q_d) \mapsto \pare{\xi_{1}(q),...,\xi_{\nc}(q)}\PERIOD
\end{equation}
We emphasize that the knowledge of  ``fast forces'' is not required, and the variables $\xi$ can be chosen arbitrarily.
If the fDOFs are not identified the method retains its approximation properties while not performing efficiently.
The fDOFs are penalized with a mass-tensor modification given by
\be\label{e:PMM}
	\pen{\MASST}(q) = \MASST + \nu^{2}\nabla_{q}\xi \,\MASSTZ\,  \nabla_{q}^T\xi \COMMA
\ee
where $\nu$ denotes the penalty intensity, and $\MASSTZ$ a ``virtual'' mass matrix associated with the fDOFs. The modification does not impact motions orthogonal to the fDOFs.
The position dependence of the mass-penalization introduces a geometric bias. This bias is corrected by introducing an effective potential
\be\label{e:Vfixnu}
V_{\M{fix},\nu}(q) = \frac{1}{2\beta} \ln \pare{ \DET (\pen{\MASST}(q)) } \COMMA
\ee
which will turn out to be a $\nu^{-1}$-perturbation of the usual Fixman corrector (see \cite{Fix78}) associated with the sub-manifold
defined by constraining the fDOFs $\xi$.
The key point is then to use an implicit representation
of the mass penalty with the aid of the extended Hamiltonian
\begin{equation}\label{e:Himmp}
\syst{
& H_{\M{IMMP}}(p,p_{z},q,z) = \frac{1}{2} p^{T} \MASST^{-1} p  + \frac{1}{2} p_{z}^{T} \MASSTZ^{-1} p_{z}
             + V(q) + V_{\M{fix},\nu}(q)\COMMA & \\
& \xi(q) = \frac{z}{\nu} \PERIOD & \TAG{\pen{C}}
}\end{equation}
The auxiliary degrees of freedom $z$ are
endowed with the ``virtual'' mass-matrix $\MASSTZ$. The constraints $(\pen{C})$ are applied in order to identify
the auxiliary variables and the fDOFs $\xi$ with a coupling intensity tuned by $\nu$. 
The typical time scale of the fDOFs is thus enforced by the penalty $\nu$.  
The system is coupled to a thermostat through a Langevin equation \eqref{e:langevin}, 
which yields a stochastically perturbed dynamics that samples the equilibrium canonical distribution. We then obtain the following desirable properties:
\begin{enumerate}
  \item The associated canonical equilibrium distribution in position is independent of the penalty 
        $\nu$.
  \item The limit of vanishing penalization ($\nu = 0$) is the original full dynamics,  enabling 
        the construction of dynamically consistent numerical schemes.
  \item The limit of infinite penalization is a standard effective constrained dynamics on the 
        ''slow'' manifold associated with stiff constraints on $\xi$.
  \item Numerical integrators can be obtained through a simple modification of standard integrators 
        for effective dynamics with constraints yielding equivalent computational complexity.
 \end{enumerate}

The dynamics associated with the IMMP Hamiltonian \eqref{e:Himmp} is detailed in \eqref{e:IMMP}, see Section~\ref{s:IMMP}.
The numerical discretization (using a leapfrog/Verlet splitting with constraints, 
usually called ``RATTLE'' for fully constrained dynamics)
is given by \eqref{e:immpscheme}. When considering sampling, the time-step error of the numerical flow can be corrected with a Metropolis step (the so-called Generalized Hybrid Monte-Carlo 
method, see references in Section~\ref{s:numerinteg}) to obtain {\it exact sampling}.
When this correction is introduced, the gradient of the Fixman potential \eqref{e:Vfixnu} need not be computed.
These numerical aspects of the method are presented in Section~\ref{s:numerinteg}. By including a 
penalty, the proposed method modifies the original Hamiltonian. However, the mass penalty can be 
also thought of as depending on the time step $\nu = \nu(\dt)$ leading to order two consistent 
schemes.

In Section~\ref{s:highoscill},  we introduce a small stiffness parameter $\epsilon$ encoding the fastest DOFs, and show that the penalty intensity can be
scaled as $\nu = \nusc/\epsilon$ in order to obtain asymptotically stable dynamics in the limit $\epsilon \to 0$. 
We prove that the dynamics converge towards the expected Markovian effective dynamics on the slow manifold. We also present analysis of the corresponding asymptotic preserving properties 
of the proposed scheme.

High-dimensional systems usually contain a large variety of timescales, and are therefore challenging test cases.
The $N$-alkane model is numerically studied in Section~\ref{s:numerexperiment} and systematically compared to Verlet scheme and constrained integration, with separate studies for dynamical and sampling issues.
In the case of butane, bond angles are penalized. 
Dynamical interpolation between exact dynamics and rigidly constrained dynamics 
is demonstrated, with the associated gain in the time step stability. 
On the other hand, exact sampling with possible gain in the mixing time for Metropolized 
sampling methods is analyzed.
Furthermore, for the $N$-alkane model with large $N$, torsion angles are penalized, a large mass-penalty 
of order $\BIGO(N)$, i.e., $\nu =\nusc N $, where $N$ is the system size is considered. For dynamics, it induces a gain in the time step 
stability region that grows with $N$, while numerical evidence is given that some macroscopic 
timescales, e.g., low frequencies  of the chain length dynamics,
remain of order $\BIGO(1)$. For sampling, it induces a similar gain for mixing time in terms of iteration steps, and measured with autocorrelation of the chain length evolution. Rigorous proofs with explicit scalings of this behavior
are provided for the case of a linear atomic chain with
quadratic (harmonic) interactions in Section~\ref{s:Nanal}, 
and consistence of the IMMP macroscopic
dynamics towards a stochastic wave equation when the re-scaled penalty 
vanishes $\nusc \to 0$ is demonstrated.

\medskip\noindent{\bf Acknowledgments:}
The research of M.R. was partially supported by the EPSRC grant GR/S70883/01 while he was visiting Mathematics Institute,
University of Warwick. The research of P.P. was partially supported by the National Science Foundation under the grant
NSF-DMS-0813893 and by the Office of Advanced Scientific Computing Research,
U.S. Department of Energy; the work was partly done at the ORNL, which is managed by UT-Battelle, LLC 
under Contract No. DE-AC05-00OR22725.

\section{Langevin processes and sampling of canonical distribution}\label{s:notations}
We consider a Hamiltonian system in the phase-space $\R^{d}\times\R^{d}$
with the Hamiltonian $H$ in the form
\be\label{e:HH}
H(p,q) = \frac{1}{2} p^{T} \MASST^{-1} p  + V(q) \COMMA
\ee
We use generic matrix notation, for instance, the Euclidean scalar product of two vectors $p_1,p_2\in\R^N$ is
denoted by $p_1^T p_2$, and the gradients of mappings from $\R^d$ to $\R^n$ with respect to standard bases
are represented by matrices
\[
 (\nabla^T_q\xi)_{ij} = (\nabla_q\xi)_{ji} = \frac{\partial \xi_i}{\partial q_j}\COMMA\; i=1,\dots,n\,,\;
                                                                                         j=1,\dots,d \PERIOD
\]

When the system is thermostatted, i.e., kept at the constant temperature, the long time distribution of the system in the
phase-space is given by the canonical equilibrium measure at the inverse temperature $\beta$
(also called the NVT distribution) given by
\be \label{e:boltzmann}
\mu(dp\,dq) = \frac{1}{Z}\EXP{ - \beta H(p,q) } dp\, dq \COMMA\;\;\; Z = \int_{\R^d\times\R^d} \EXP{ - \beta H(p,q) } dp\,dq\COMMA
\ee
with the normalization constant $Z<\infty$. The standard dynamics used to model thermostatted systems
are given by Langevin processes.
\begin{Def}[Langevin process]\label{d:langevin} A Langevin process at the inverse temperature $\beta$ with the Hamiltonian
  $H(q,p)$, $(p,q)\in\R^d\times \R^d$, the $d\times d$ dissipation matrix $\gamma$, and the fluctuation matrix $\sigma$
  is given by the stochastic differential equations
  \be\label{e:langevin}
  \syst{
    &\dot{q} = \Dp H  &\\
    &\dot{p} = -\Dq H -\gamma \dot{q} + \sigma \dot{W} \COMMA &
    }
   \ee
  where $\dot{W}$
  is a standard white noise (Wiener process), and
  $\sigma\in\R^d\times\R^d$
  satisfies the fluctuation-dissipation identity
  \[
  \sigma \sigma^{T} =\frac{2}{\beta} \gamma \PERIOD
  \]
For any $\gamma$, the process is reversible with respect to the stationary canonical distribution \VIZ{e:boltzmann}.
Furthermore, if $\gamma$ is strictly positive definite, the process  is ergodic.
\end{Def}

Throughout the paper, stochastic integrands have finite variation
thus the stochastic integration (e.g., It\^o or Stratonovitch) need not be specified. Furthermore, the usual global
Lipschitz conditions  (see \cite{Osk92}) on $H$ and $\xi$ are assumed, ensuring well-posedness of the considered stochastic 
differential equations. The analysis presented in the paper can be generalized to a position dependent
dissipation matrix $\gamma = \gamma(q)$.

The mapping $\xi:\R^d\to\R^n$, defines
$\nc \leq d$ degrees of freedom, given by smooth  functions
taking values in a neighborhood of $0$. We assume that the mapping $\xi$ is regular
(i.e., with a non-degenerate Jacobian) in an open $\delta$-neighborhood
$\mathcal{O}_{\delta}=\set{q\SEP \norm{\xi(q)} < \delta} $ of $\xi^{-1}(0)$, hence defining a smooth sub-manifold of
$\mathbb{R}^{d}$ denoted $\SMANZ{z}=\xi^{-1}(z)$ for $z$ in a neighborhood of the origin.
The dependence of the potential $V$ with respect to the degrees of freedom $\xi$ is expected to be ``stiff'' in the second variable. In Section \ref{s:highoscill} we will introduce the stiffness parameter
$\epsilon$. In that section we shall assume that such parameter dependence can be 
explicitly identified, and that the potential energy
$V$ can be written in the form
\begin{equation}\label{Voscill}
V(q)= U(q,\frac{\xi(q)}{\epsilon}) \COMMA
\end{equation}
where the function $U:\R^d\times\R^n\to\R$ satisfies the coercivity condition
$\lim_{z \to +\infty} U(q,z) = +\infty$. The fast degrees of freedom $\xi$
of states at a given energy then remain in a closed neighborhood of the origin as the
stiffness parameter $\epsilon \to 0$. In this limit the system is confined to the sub-manifold $\SMAN$ which is usually called the ``slow manifold''.

\section{The implicit mass-matrix penalization method}\label{s:IMMP}
In this section we focus on properties of the IMMP method. The multiscale
structure of the potential $V$ need not be known in order to apply the method.
Thus, in this section, we consider the potential $V$ in the form where we do not impose the
structural assumption \VIZ{Voscill} on the potential function $V:\R^d\to\R$.

\subsection{Description of the method}
The new, penalized mass-matrix of the system is the position dependent tensor defined in \eqref{e:PMM}. The associated modified impulses are denoted
\be\label{e:ppen}
	\pen{p} = \pen{\MASST}(q) \MASST^{-1} p \PERIOD
\ee
When $\nu$ becomes large, the velocities are bound to remain tangent to the manifolds $\set{q|\xi(q) = z }$,
and orthogonal motions are arbitrarily
slown down. Conversely, when $\nu=0$, one recovers the original highly oscillatory system.
Since the modification in $\pen{\MASST}$ depends on the position $q$,  new geometry is introduced and an additional correction \eqref{e:Vfixnu} in the potential energy is required in order to preserve original statistics in the position variable. This correction is in fact close to the standard Fixman corrector  for $\nu$ large (see \eqref{e:Vfix}).
Defining $G(q)$ as the $\nc \times \nc$ Gram matrix associated with the fast degrees of freedom
\begin{equation}\label{e:Gram}
G(q) =   \nabla^{T}_{q} \xi \,\MASST^{-1}  \Dq \xi \COMMA
\end{equation}
one has the following property of the correcting potential.
\begin{Pro} Up to an additive constant, we have
\be\label{e:penV}
V_{\M{fix},\nu}(q) =  \frac{1}{2\beta} \ln \DET \pare{ G(q) +   \frac{1}{\nu^{2} } \MASSTZ^{-1} } \COMMA
\ee
and thus (up to additive constants)
\[
\lim_{\nu \to +\infty} V_{\M{fix},\nu} = V_{\M{fix}} = \frac{1}{2\beta} \ln \DET \pare{ G(q) } \COMMA\;\;\;\;\mbox{and }\;\;\;\;\;
\lim_{\nu \to 0} V_{\M{fix},\nu} = 0 \PERIOD
\]
\end{Pro}
\begin{proof}
Using the  identity for a non-diagonal matrix $J$ of dimension $n_{1} \times n_{2}$:
\[
  \DET( \ID_{n_1}+ JJ^{T}   ) = \DET( \ID_{n_2}+ J^{T}J   ) \COMMA
\]
one observes
\[
  \DET(\pen{\MASST}) = \DET(\MASST)\, \DET(\nu ^{2}\MASSTZ)\, \DET(G+\frac{1}{\nu^2}\MASSTZ^{-1})
\]
from which the expression for the corrected Fixman potential follows.
\end{proof}

The associated modified Hamiltonian is then given by
\be\label{e:Hnu}
	\pen{H}(\pen{p},q) = \frac{1}{2} \pen{p}^{T} {\pen{\MASST}}^{-1} \pen{p}    
        + V(q) + V_{\M{fix},\nu}(q)  \COMMA
\ee
and $H_{0}=H$ is the original Hamiltonian \eqref{e:H}.

Sampling such a system can be done using the standard Langevin stochastic perturbation
as detailed in Definition~\ref{d:langevin}.
However, the direct discretization of the equation of motion given by $\pen{H}$ (e.g., by an explicit scheme)
is bound to be unstable from non-linear instabilities when the fast degrees of freedom are not affine functions.
In order to construct stable schemes one may rather use an implicit formulation of the Hamiltonian \eqref{e:Hnu},
in conjunction with a solver which enforces the constraints.
To obtain such a formulation we extend the state space with $\nc$ new variables $(z_{1},..,z_{\nc})$,
and associated moments $(p_{z_{1}},..,p_{z_{\nc}})$. The auxiliary mass-matrix for the
new degrees of freedom is then given by $\MASSTZ$.
The new extended Hamiltonian of the system $H_{\M{IMMP}}$, defined by \eqref{e:Himmp}, is now defined in 
$\mathbb{R}^{d+\nc}\times\R^{d+\nc}$, where $\nc$ position constraints denoted by $(\pen{C})$ are included.
This construction implies $\nc$ hidden constraints on momenta. The equivalence of the two Hamiltonians
\eqref{e:Hnu} and  \eqref{e:Himmp} formulations is stated as a simple separate lemma.
\begin{Lem}
The equations of motion associated with the penalized mass-matrix Hamiltonian \eqref{e:Hnu} or the extended Hamiltonian
with constraints \eqref{e:Himmp} are identical.
\end{Lem}
\begin{proof}
The Lagrangian associated with $H_{\M{IMMP}}$ is given by
\[
   L_{\M{IMMP}}(\dot{q},\dot{z},q,z) = \frac{1}{2} \dot{q}^{T} \MASST \dot{q}  +\frac{1}{2} 
   \dot{z}^{T} \MASSTZ \dot{z} - V(q) -  V_{\M{fix},\nu}(q) \COMMA
\]
and includes hidden constraints on velocities $\dot{z} = \nu \Dq^T \xi \,\dot{q}$
implied by the constraints $(\pen{C})$ on position variables.
Replacing $\dot{z}$ and $z$ in $L_{\M{IMMP}}$ by their expressions as functions of $\dot{q}$ and $q$,
one obtains the Lagrangian associated with $\pen{H}$.
\end{proof}

The stochastically perturbed equations of motion of the Langevin type associated with \eqref{e:Himmp} define
the dynamics with implicit mass-matrix penalization.
\begin{Def}[IMMP]\label{d:IMMP} The implicit Langevin process associated with Hamiltonian $H_{\M{IMMP}}$
and constraints $(\pen{C})$ is defined by the following equations of motion
\begin{equation} \label{e:IMMP}\syst{
    &\dot{q} =  \MASST^{-1} p  &  \\
    &\dot{z} =  \MASSTZ^{-1} p_{z} &  \\
    &\dot{p} = -\Dq V(q) - \Dq V_{\M{fix},\nu}(q) -\gamma \dot{q} + \sigma \dot{W} - \Dq \xi 
                \,\dot{\lambda} & \\
    &{\dot{p}_{z}} = -\gamma_{z} \dot{z} + \sigma_{z} \dot{W}_{z} +  \frac{\dot{\lambda}}{\nu} & \\
    &\xi(q) = \frac{z}{\nu}\COMMA & \TAG{\pen{C}}
  }
\end{equation}
The process $\dot{W}$ (resp. $\dot{W}_{z}$ ) is a standard multi-dimensional white noise, $\gamma$ 
(resp. $\gamma_{z}$) a $d\times d$ (resp. $\nc\times \nc$) non-negative symmetric 
dissipation matrix, $\sigma$ (resp. $\sigma_{z}$) is the fluctuation matrix satisfying 
$\sigma \sigma^{T} =\frac{2}{\beta} \gamma$ 
(resp. $\sigma_{z} \sigma_{z}^{T} =\frac{2}{\beta} \gamma_{z}$).
The processes $\lambda \in \mathbb{R}^{\nc}$ are Lagrange multipliers associated 
with the constraints $(\pen{C})$ and adapted with the white noise.
\end{Def}

This process is naturally equivalent to the explicit mass-penalized Langevin process in $\mathbb{R}^{d}\times\R^d$
associated with $\pen{\MASST}$. 
Moreover, when the penalization vanishes ($\nu \to 0$), the evolution law of
the process $\{p_t,q_t\}_{t\geq 0}$ or $\{({p_{\nu}})_t,q_{t}\}_{t\geq 0}$ converges towards the original dynamics.
\begin{Pro}\label{p:IMMPeq} The stochastic process with constraints \eqref{e:IMMP} is well-posed and equivalent
to the Langevin diffusion in $\mathbb{R}^{d}\times\R^d$ (see Definition~\ref{d:langevin}),
with the mass-penalized Hamiltonian $\pen{H}$ \eqref{e:Hnu}, and the dissipation matrix given by
\[
\pen{\gamma}(q) = \gamma + \nu^{2} \Dq\xi \,\gamma_{z} \nabla^{T}_{q}\xi \PERIOD %
\]
Furthermore, the process is reversible and ergodic with respect to the canonical distribution \eqref{e:penboltz} (with marginal in position variables given by the original potential, i.e., up to the normalization,
$\EXP{ - \beta V(q) } dq$) .
\end{Pro}
\begin{proof}
Imposing the  constraints implies $ \nabla_q^T \xi\, \MASST^{-1} p = \frac{1}{\nu} \MASSTZ^{-1} p_{z}$.
Thus by the definition of $\pen{p}$ we have
\[
 \pen{p} = p + \nu \Dq\xi \, p_z \PERIOD
\]
 Since the position process $\{q_{t}\}_{t\geq 0}$ is of finite variation, a short computation shows that for each coordinate $i = 1,..,d$
\be\label{e:step}
\pen{\dot{p}}^i = {\dot{p}^i} + \nu \partial_{q_i} \xi \, \dot{p}_{z} + \nu^{2} \dot{q}^T   \Dq \! (\partial_{q_i}\xi) \, p_z \PERIOD
\ee
Furthermore,
\[
  -\partial_{q_i} \pare{\frac{1}{2}  \pen{p}^T \pen{\MASST}^{-1} \pen{p}}= \partial_{q_i}\pare{\frac{1}{2} \dot{q}\pen{\MASST} \dot{q}}
  = \nu^2 \dot{q}^T  \Dq(\!\partial_{q_i} \xi) \MASSTZ\nabla_q^T \xi\, \dot{q} \COMMA
\]
and thus
\[
{\pen{\dot{p}}} = \dot{p} + \nu \Dq \xi \,{\dot{p}_{z}} - \Dq \pare{\frac{1}{2}  \pen{p}^T \pen{\MASST}^{-1} \pen{p}} \PERIOD
\]
Substituting the expressions for  $\dot{p}$ and ${\dot{p}_{z}}$ from \eqref{e:IMMP} into \eqref{e:step} we obtain
\begin{equation}\label{e:IMMPeq}
{\pen{\dot p}}  =  -\frac{1}{2} \pen{p}^{T} \Dq \pen{\MASST}^{-1} \pen{p} - \Dq V(q)  - \Dq V_{\M{fix},\nu}(q)
\, -\gamma \dot{q} - \nu \Dq \xi\, \gamma_{z} \dot{z} + \sigma \dot{W} + \nu \Dq \xi\, \sigma_{z} \dot{W}_{z} \COMMA
\end{equation}
which yields the result.
\end{proof}
\subsection{Exact sampling in position variables}
By construction, statistics of positions $q$ of the mass penalized Hamiltonian are independent of the penalization,
leading to the \emph{exact canonical statistics} in  position variables.
\bpro[Exact statistics]
The canonical distribution associated with the mass-penalized Hamiltonian \eqref{e:Hnu} is given by
\be\label{e:penboltz}
\pen{\mu}(d\pen{p} \, dq) = \frac{1}{\pen{Z}} \EXP{ - \beta \pen{H}(\pen{p},q) } d\pen{p}\,dq \PERIOD
\ee
Its marginal probability distribution in $q$ is
\[
 \frac{1}{\pen{Z}} \int  \EXP{ - \beta \pen{H}(\pen{p},q) } d\pen{p} = 
     \frac{\EXP{ - \beta V(q) } dq}{\int \EXP{ - \beta V(q) }\,dq}
\]
which is the original canonical distribution \eqref{e:boltzmann} in the position variables,
and is independent of the mass penalization parameter $\nu$.
\epro
\begin{proof}
The normalization of Gaussian integrals in the $\pen{p}$ variables yields
\[
\int \EXP{ - \beta \frac{1}{2} \pen{p}^{T} {\pen{\MASST}}^{-1} \pen{p} } d \pen{p} =
 \left(\frac{2 \pi}{\beta}\right)^{d/2}\!\sqrt{\DET(\pen{\MASST})} \COMMA
\]
which is cancelled out  by the Fixman corrector $V_{\M{fix},\nu}$ and the result follows.
\end{proof}
\subsection{Interpolation between exact and constrained dynamics}
In this section, the IMMP dynamics is shown to be an interpolation between exact dynamics ($\nu=0$), and 
constrained dynamics ($\nu=+\infty$).

\begin{Pro}[Small penalty]\label{p:vanishpen}
When $\nu \to 0$ the evolution law of the processes $\{p_t,q_t\}_{t\geq 0}$ or
$\{(\pen{p})_t,q_t\}_{t\geq 0}$ defined by the implicit equations \eqref{e:IMMP} converges
(in the sense of probability distributions on continuous paths endowed with the uniform convergence)
towards the process solving the original Langevin dynamics \eqref{e:langevin}.
\end{Pro}
 \begin{proof}
The stochastic differential equation defined by $\dot{q} = \pen{M}^{-1} \pen{p}$
and \eqref{e:IMMPeq} has smooth coefficients which depend on $\nu$ in a continuous fashion ($\nu \mapsto \pen{M}$
and $\nu  \mapsto V_{\M{fix},\nu}$ are continuous).
Standard results on weak convergence (\cite{EthKur86}) of stochastic processes imply the result as stated.
\end{proof}

When the mass penalty tends to infinity, the IMMP process converges to a constrained process on the manifold $\SMANZ{z_{t=0}}=\set{ q \SEP \xi(q)=z_{t=0}}$. 
\bpro[Large penalty]\label{p:largepen}
        Consider a family of initial conditions indexed by $\nu$ and satisfying
	\[
	\sup_{\nu} \abs{ \nu\,(\xi(q_{t=0}) - z_{t=0}) } < +\infty \COMMA
	\]
        and assume that the Gram matrix $G$ is invertible in a neighborhood of $ \SMANZ{z_{t=0}}$.
	Then when $\nu \to +\infty$ the IMMP Langevin stochastic process \eqref{e:IMMP} converges in distribution
        towards the decoupled limiting processes with constraints
	\begin{equation}\label{e:largepen}
        \syst{
    	  &\dot{q} =  \MASST^{-1} p \COMMA &  \\
    	  &\dot{p} = -\Dq V  - \Dq V_{\M{fix}} -\gamma \dot{q} + \sigma \dot{W} - \Dq\xi \dot{\lambda}\COMMA  & \\
    	  &\xi(q)  = z_{t=0}    \COMMA                                        &  \TAGG{C} \\
    	  &\dot{z} =  \MASSTZ^{-1} p_{z} \COMMA &  \\
    	  &{\dot{p}_{z}} = -\gamma_{z} \dot{z} + \sigma_{z} \dot{W}_{z} \PERIOD &
	}
        \end{equation}
	where $\{\lambda_t\}_{t\geq 0}$ are adapted stochastic processes defining
        the Lagrange multipliers associated with the constraints $(C)$.

	Furthermore, the process $\{q_t,p_t\}_{t\geq 0}$ defines an effective dynamics with constraints
        (see also Definition~\ref{d:consteff}) on the sub-manifold $\SMANZ{z_{t=0}}$. It is reversible with respect to 
its stationary canonical distribution given, up to the normalization, by the ``stiff'' Boltzmann distribution
        \[
	\EXP{-\beta (H(p,q)+V_{\M{fix}}(q))} \sigma_{T^{*}\SMANZ{z_{t=0}}}(dp\,dq)
	\]
	with the $q$-marginal $\EXP{- \beta V(q)}\,\delta_{\xi(q) = 0}(dq)$.
	When $\gamma$ and $\gamma_{z}$ are strictly positive definite the process is ergodic.
\end{Pro}
\begin{proof}
By a simple translation, it is sufficient to show the proposition for $z_{t=0} = 0$. Satisfying the constraint $(\pen{C})$ 
in \eqref{e:IMMP} implies a hidden constraint in the momentum space,
$\nabla_q \xi \MASST^{-1} p$ $= \frac{1}{\nu} \MASSTZ^{-1} p_{z}$. Differentiating this expression with respect to time
and replacing the result in \eqref{e:IMMP} yields an explicit formula for the Lagrange multipliers
\be\label{e:lagecd}
  \dot{\lambda} = (G+\frac{1}{\nu^2}\MASSTZ^{-1})^{-1}\left[\HESS(\xi) \pare{\MASST^{-1} p,\MASST^{-1} p} +
                                                             \Dq \xi \MASST^{-1} f_{q}
                                                           - \frac{1}{\nu}\MASSTZ^{-1} f_{z}\right] \COMMA
\ee
with forces $(f_q,f_z)$
\begin{eqnarray*}
f_{q} &=&  -\Dq V - \Dq V_{\M{fix},\nu}  -\gamma \MASST^{-1} p + \sigma \dot{W}\COMMA \\
f_z &=&   - \gamma_{z} \MASSTZ^{-1} p_{z} + \sigma_{z} \dot{W}_{z}\COMMA
\end{eqnarray*}
and the Hessian $\HESS(\xi)$ of the mapping $\xi$ acting on the velocities $\MASST^{-1}p$.
This calculation shows that \eqref{e:IMMP} is in fact a standard stochastic differential equation with smooth coefficients,
and thus has a unique strong solution. The coefficients of these stochastic differential equations are
continuous in the limit $\tfrac{1}{\nu} \to 0$, at least in a $\delta$-neighborhood of $\SMAN$ in
which $G$ is invertible.
The formally computed  limiting process is given by \eqref{e:largepen} with the Lagrange multipliers solving
 \[
 \dot{\lambda} = G^{-1}\pare{\HESS(\xi) \pare{\MASST^{-1} p,\MASST^{-1} p} + \Dq^T \xi \, \MASST^{-1} f_{q}} \PERIOD
 \]
By construction, this limiting process satisfies the constraint $\xi(q)=0$.
Its coefficients are Lipschitz and the process is well-posed.
As a result of those properties, the rigorous proof of weak convergence follows
classical arguments, see \cite{EthKur86}, that are divided into three steps
\begin{description}
\item[{\rm (i)}] We truncate the process \eqref{e:IMMP} to a compact neighborhood of $\mathcal{M}_{0}$.
\item[{\rm (ii)}] The continuity of the Markov generator with respect to $\tfrac{1}{\nu}$ implies tightness for 
the associated $\tfrac{1}{\nu}$-sequences of truncated processes with the limit being uniquely defined by \eqref{e:largepen}.
\item[{\rm (iii)}] The limiting process remains on  $\SMAN$, which implies weak convergence of the sequence without truncation.
\end{description}
The process \eqref{e:largepen} is thus a Langevin process with constraints, exhibiting reversibility properties
with respect to the associated Boltzmann canonical measure and is ergodic when $\gamma$ is strictly positive definite
(see the summary in Appendix~\ref{s:langevinapp}).
Note that the $q$-marginal is geometrically corrected by the Fixman potential term.
\end{proof}

We conclude this section by discussing some consequences for numerical computations.
\begin{Rem}
{\rm
Proposition~\ref{p:vanishpen}  and ~\ref{p:largepen} imply that the IMMP 
     scheme is a tunable
     interpolation between the exact stochastic dynamics \eqref{e:langevin}, and the stochastic 
     dynamics with constraints \eqref{e:largepen}. If one prefers to interpolate with rigidly 
     constrained dynamics (i.e., without the Fixman correction, see also 
     Section~\ref{s:highoscill} for a detailed discussion on ``stiff'', as opposed 
     to ``rigid'', constrained dynamics), one removes the Fixman correction in the force 
     evaluation.
}
\end{Rem}

\section{Numerical integration}\label{s:numerinteg}
The key ingredient for achieving efficient numerical simulation is to use an integrator
that enforces the constraints associated with the implicit formulation of the mass penalized dynamics \eqref{e:IMMP}. 
The implicit structure of \eqref{e:IMMP} leads to numerical schemes that are potentially asymptotically stable in stiff 
cases (Section~\ref{s:highoscill}). On the other hand, when the penalization $\nu$ vanishes with the time-step
the scheme becomes consistent with respect to
the original exact dynamics \eqref{e:langevin}. One may then consider the mass-penalization introduced here as a special method of
pre-conditioning for a stiff ODE system with an ``implicit'', in the time evolution sense, structure.
Here, the ``implicit'' structure amounts to solving the imposed constraints $\xi(q) = z/\nu$ in \VIZ{e:IMMP}.

It lies outside the scope of this paper to review standard numerical methods for 
constrained mechanical systems,
we refer to \cite{HaiLub02} as a classical textbook, and to the series 
(\cite{GunBer77, RycCic77, CicVan06, CicKapVan05, CicLelVan06})
as a sample of works on practical developments of numerical methods.
The IMMP method is presented  with the classical leapfrog/Verlet scheme that enforces constraints,
usually called RATTLE in \eqref{e:immpscheme}. It can be implemented by a simple modification 
of standard schemes constraining fDOFs. The scheme is second order, reversible and symplectic. 
This choice is largely a presentation matter, for practical purposes one can refer to one's
favorite numerical integrator for Hamiltonian systems with or without stochastic perturbations. 

For accurate sampling of the equilibrium distribution, one can also add a Metropolis
acceptance, rejection time-step corrector 
at each time step of the deterministic integrator. If the underlying integrator is reversible and preserves 
the phase-space measure, this extension leads to a scheme which {\it exactly} preserves canonical distributions. 
The Metropolis correction is used in Hybrid Monte-Carlo (HMC) methods,
which are sampling algorithms relying on the underlying dynamics of the system to generate moves in 
the configuration space, which are accepted or rejected according to the Metropolis rule. The 
Metropolis acceptance/rejection step can be used at each integration
time step of the Langevin process, usually referred to as Generalized Hybrid Monte Carlo
(GHMC) introduced in \cite{Horowitz91}. 
However, the necessity of the momentum flip when a rejection occurs destroys the dynamical 
features of the Langevin process when the rate of rejection does not vanish, and makes the 
latter rather behave similarly to an overdamped dynamics. Note also that for usual schemes 
integrating the Hamiltonian dynamics, the average acceptance ratio often decreases when the 
dimension of the system increases (\cite{IzaHam04}).  Many improvements and modifications 
(\cite{Creutz89, Horowitz91, Beccaria94, Kennedy98}) of the HMC algorithm  have been developed 
since its introduction in \cite{Duane85} for simulations applied to quantum statistical field theories. 
Subsequently it has been also employed to a wide range of simulations in macromolecular
systems, e.g., \cite{ClaBakStiBra94,SchFis99}.

We implement numerical discretization of the Langevin process with constraints \eqref{e:IMMP} obtained
by splitting the Hamiltonian part and the Gaussian fluctuation/dissipation perturbation. Note that 
by using a Metropolis acceptance/rejection rule (HMC), or simply by weighting statistical averages,
the forces associated with the Fixman corrector need not be computed when one is interested in 
sampling only.
\subsection{The IMMP integrator}
We recall that we consider the IMMP dynamics \eqref{e:IMMP}, which consists of the following elements
\begin{enumerate}
\item the Hamiltonian $H_{\M{IMMP}}$ defined in \eqref{e:Himmp}, which defines the deterministic dynamics of the IMMP.
\item the dissipation matrix $\mathrm{diag}\,(\gamma,\gamma_z)$,
and the inverse temperature $\beta$, which defines the features of the stochastic thermostat.
\end{enumerate}

\begin{scheme}[Dynamical integrator]\label{d:scheme}
\bi
\item[Step 1:] Integrate the Hamiltonian part with:
\begin{eqnarray} \label{e:immpscheme}
		&&\syst{
		&p_{n+1/2} = p_{n} -  \frac{\dt}{2} (\Dq V+\Dq V_{\M{fix},\nu} )(q_{n}) - \Dq \xi(q_{n})  \lambda_{n+1/2}&\\
		&p_{n+1/2}^z = p_{n}^z +\frac{1}{\nu}\lambda_{n+1/2} &
		} \nonumber \\
		&&\syst{
		&q _{n+1} =  q_{n}+ \dt \MASST^{-1} p_{n+1/2} &\\
		& z_{n+1} =  z_{n}+ \dt \MASSTZ^{-1} p_{n+1/2}^z &
		} \nonumber \\
		&&\hspace{.5cm} \xi(q_{n+1}) = \frac{z_{n+1}}{\nu} \qquad \hspace{3cm} \TAG{C_{1/2}} \\
		&&\syst{
		&p_{n+1}   = p_{n+1/2} -  \frac{\dt}{2} (\Dq V+\Dq V_{\M{fix},\nu})(q_{n+1})  - \Dq \xi(q_{n+1})  \lambda_{n+1}&\\
		&p_{n+1}^z = p_{n+1}^z +\frac{1}{\nu}\lambda_{n+1}&
		} \nonumber \\
		&&\hspace{.5cm} \Dq^T \xi(q_{n+1})\,\MASST^{-1}p_{n+1} = \frac{1}{\nu} \MASSTZ^{-1}p^z_{n+1}\qquad \hspace{3cm}
                 \TAG{C_{1}}   \PERIOD \nonumber
\end{eqnarray}
\item[Step 2:] Integrate if necessary the Gaussian fluctuation/dissipation part with a mid-point 
Euler scheme with constraints (see Appendix~\ref{s:exactflucdiss}).
\ei
\end{scheme}
Here again, one can remove the Fixman correction forces $\Dq V_{\M{fix},\nu}$ in force 
evaluation in order to interpolate with usual rigid constraints dynamics.

Note a useful variant of the integrator, which occurs when the potential dependence with 
respect to the penalized variables is known explicitly
$V(q) = U(q,\xi(q))$.
In such a case, the expressions \VIZ{e:immpscheme} is replaced by
\begin{eqnarray} \label{e:immpscheme2}
		&&\syst{
		&p_{n+1/2} = p_{n} -  \frac{\dt}{2} (\Done U+\Dq V_{\M{fix},\nu} )(q_{n}) - \Dq \xi(q_{n})  \lambda_{n+1/2}&\\
		&p_{n+1/2}^z = p_{n}^z -  \frac{\dt}{2\nu}\Dtwo U(q_n) + \frac{1}{\nu}\lambda_{n+1/2} &
		} \nonumber \\
		&&\syst{
		&q _{n+1} =  q_{n} + \dt \MASST^{-1} p_{n+1/2} &\\
		& z_{n+1} =  z_{n}+ \dt \MASSTZ^{-1} p_{n+1/2}^z &
		} \nonumber \\
		&&\hspace{.5cm} \xi(q_{n+1}) = \frac{z_{n+1}}{\nu} \qquad \hspace{3cm} \TAG{C_{1/2}} \\
		&&\syst{
		&p_{n+1}   = p_{n+1/2} -  \frac{\dt}{2} (\Done U+\Dq V_{\M{fix},\nu})(q_{n+1})  - \Dq \xi(q_{n+1})  \lambda_{n+1}&\\
		&p_{n+1}^z = p_{n+1}^z -  \frac{\dt}{2\nu} \Dtwo U(q_{n+1}) +\frac{1}{\nu}\lambda_{n+1}&
		} \nonumber \\
		&&\hspace{.5cm} \Dq^T \xi(q_{n+1})\,\MASST^{-1}p_{n+1} = \frac{1}{\nu} \MASSTZ^{-1}p^z_{n+1}\qquad \hspace{3cm}
                 \TAG{C_{1}}   \PERIOD \nonumber
\end{eqnarray}
where in the above, $\Done$ and $\Dtwo$ denote the derivatives with respect two the first and the 
second variables, respectively. This variant consists in applying a part of the force to the 
auxiliary variables instead of directly to the system. For penalized degrees of freedom that move far 
from their equilibrium value, e.g., torsion angles, see Section~\ref{s:numerexperiment}, this may 
increase the stability of the algorithm enforcing the constraints. In the same way for slow/fast 
systems (see Section~\ref{s:infstiff}), the scheme \VIZ{e:immpscheme2} will lead to 
asymptotic stability in the large stiffness limit.

\subsection{Exact sampling and Monte-Carlo scheme}
To obtain exact sampling by correcting time step errors and, if necessary, the geometric bias, a Metropolis acceptance/rejection is added. A domain of phase-space $D_{\dt} \subset \R^{2d} \times \R^{2n}$ where the integrator \eqref{e:immpscheme} with constraints has a unique solution is 
also considered. In practice this set is simply the set of configurations for which the algorithm used to enforce constraints in \eqref{e:immpscheme} (e.g., a Newton algorithm) converges in a given number of steps. In practice, the (equilibrium) probability for the system of lying outside $D_{\dt}$ goes to zero exponentially fast with the time step, and the latter is taken such that this probability is negligible.
\begin{scheme}[Leapfrog/Verlet algorithm with Metropolis correction for Langevin IMMP \eqref{e:IMMP}]\label{d:schemehmc}
\bi

\item[Step 1:] Compute $(q_{n+1},z_{n+1},p_{n+1},p^{z}_{n+1})$ with the integrator \eqref{e:immpscheme}, and set 
$$ 
\Delta H_{n+1} = H_{\M{IMMP}}(q_{n+1},z_{n+1},p_{n+1},p^{z}_{n+1}) - H_{\M{IMMP}}(q_{n},z_{n},p_{n},p^{z}_{n})\PERIOD
$$
If $(q_{n+1},z_{n+1},p_{n+1},p^{z}_{n+1})$ does not belong
to $D_{\dt}$, set $\Delta H_{n+1} = + \infty$.  
\item[Step 2:] Accept the step with the probability $\min(1,{e \rm }^{-\beta \Delta H_{n+1} })$, otherwise reject, flip momenta, and set
\[
(q_{n+1},z_{n+1},p_{n+1},p^{z}_{n+1})= (q_{n},z_{n},-p_{n},-p^{z}_{n}) \PERIOD
\]

\item[Step 3:] Integrate the Gaussian fluctuation/dissipation part with a mid-point Euler scheme (for details see Appendix~\ref{s:exactflucdiss}).
\ei
\end{scheme}
\begin{Rem}
{\rm
A useful variant, when using HMC strategies, consists in modifying the Hamiltonian $H_{\M{IMMP}}$ \eqref{e:Himmp} in the integrator (importance sampling) by neglecting the Fixman corrector
\begin{equation}\label{e:Hnum}
  \syst{
  &\tilde{H}_{\M{IMMP}}(p,p_{z},q,z) =  \frac{1}{2} p^{T} \MASST^{-1} p  + \frac{1}{2} p_{z}^{T} \MASSTZ^{-1} p_{z} + \tilde{V}(q)  & \\
  &\xi(q) = \frac{z}{\nu}         & \TAG{\pen{C}},
  }
\end{equation}
where $\tilde{V}$ is a potential that may be chosen arbitrarily (typically $\tilde{V} = V$). Indeed, only the underlying phase space structure, 
which does not depend on $\tilde{V}$, is necessary in HMC methods. The correct potential $V + V_{\M{fix},\nu}$ has to be used 
in the Metropolis step only, or alternatively by weighting ensemble averages, in order to ensure exact canonical sampling. 
Thus potentially costly evaluations of the gradient of the Fixman corrector $V_{\M{fix},\nu}$ can be avoided.
This numerically constructed Markov chain preserves the canonical distribution.
}
\end{Rem}
\bpro[Exact sampling]
The Monte-Carlo algorithm generated from \eqref{e:IMMP} as
described in Scheme~\ref{d:schemehmc} generates a Markov chain that leaves invariant the canonical
distribution \eqref{e:penboltz} (conditioned in the constraints stability domain $D_{\dt}$ of \eqref{e:immpscheme}). The marginal in position variables of the distribution \eqref{e:penboltz} is the original distribution 
$\EXP{ - \beta V(q) }dq$  which is independent of the mass-penalization $\nu$.
\epro
\begin{proof}
The statement follows from reversibility and measure preserving properties of Verlet schemes (see \cite{HaiLub02}),  from the Hybrid Monte Carlo rule (\cite{DuaneKennedy87,IzaHam04}),
and from the construction of the mass-penalized Hamiltonian (Proposition~\ref{p:IMMPeq}).
\end{proof}

\subsection{Interpolation between exact and constrained dynamics}
First we study the limit towards the leapfrog/Verlet scheme.
\bpro[Small penalty limit] \label{p:shemesmallnu}
Assume that the IMMP integrator \eqref{e:immpscheme} or \eqref{e:immpscheme2} is locally well-defined, and that $\nu \to 0$. Then the latter two (computed with or without the Fixman correcting forces $\Dq V_{\M{fix},\nu}$) converge at order $\BIGO(\nu^2)$ towards the leapfrog/Verlet scheme for the exact dynamics.
\epro
\begin{proof}
Consider the auxiliary variables in the scaling $(\bar{z}, \bar{p}^z)=(\frac{z}{\nu}, \frac{p^z}{\nu})$, and the following shift of Lagrange multipliers in \eqref{e:immpscheme2}, or equivalently in \eqref{e:immpscheme}
\[
\lambda \to -\nu^2\lambda +\frac{\dt}{2} \Dtwo U,
\]
then the scheme \eqref{e:immpscheme2} becomes
 \begin{eqnarray} \label{e:immpschemenusmall}
		&&\syst{
		&p_{n+1/2} = p_{n} -  \frac{\dt}{2} (\Dq U+\Dq V_{\M{fix},\nu} )(q_{n}) + \nu^2\Dq \xi(q_{n})  \lambda_{n+1/2}&\\
		&\bar{p}_{n+1/2}^z = \bar{p}_{n}^z - \lambda_{n+1/2} &
		} \nonumber\\
		&&\syst{
		&q _{n+1} =  q_{n} + \dt \MASST^{-1} p_{n+1/2} &\\
		& \bar{z}_{n+1} =  \bar{z}_{n}+ \dt \MASSTZ^{-1} \bar{p}_{n+1/2}^z &
		}\nonumber \\
		&&\hspace{.5cm} \xi(q_{n+1}) = \bar{z}_{n+1}\qquad \hspace{3cm} \TAG{C_{1/2}} \\
		&&\syst{
		&p_{n+1}   = p_{n+1/2} -  \frac{\dt}{2} (\Done U+\Dq V_{\M{fix},\nu})(q_{n+1})  + \nu^2\Dq \xi(q_{n+1})  \lambda_{n+1}&\\
		&\bar{p}_{n+1}^z = \bar{p}_{n+1}^z  - \lambda_{n+1}&
		}\nonumber \\
		&&\hspace{.5cm} \Dq^T \xi(q_{n+1})\,\MASST^{-1}p_{n+1} = \MASSTZ^{-1}\bar{p}^z_{n+1}\qquad \hspace{3cm}
                 \TAG{C_{1}}   \COMMA \nonumber
\end{eqnarray}
which converges at order $\BIGO(\nu^2)$ to a decoupled scheme where the $(q,p)$ variables evolve according to the usual leapfrog/Verlet scheme, and the auxiliary variables $(\bar{z}, \bar{p}^z)$ are enforced by the constraints $\xi(q) = \bar{z}$. 
Now from \eqref{e:PMM}, $V_{\M{fix},\nu}$ is of order $\BIGO(\nu^2)$ which completes the proof.
\end{proof}

One can now construct schemes consistent with respect to the exact dynamics
by letting the penalty $\nu =\nusc \dt^{k}$ go to zero with
the time step, for some $k >0$.
Indeed, Proposition~\ref{p:vanishpen} shows that the mass-penalized dynamics \eqref{e:IMMP} converges towards the exact original dynamics for $\nu=0$ at order $\BIGO(\nu^2)$. 
Consequently most of the usual numerical schemes will be consistent at their own approximation order
but bounded above by $2k$. 
The order of convergence refers to the maximal integer $k$ such that the convergence of trajectories
with respect to the uniform norm occurs at the rate $\BIGO(\dt^{k})$.
Neglecting the order of the fluctuation/dissipation part (Step 2 of Scheme~\ref{d:scheme}), we deduce the following consistence property.

\bpro[Time-step consistency with exact dynamics]
Assume that the integrator \eqref{e:immpscheme} or \eqref{e:immpscheme2} are locally well-defined, and that $\nu =  \nusc \dt$. 
Then the IMMP numerical scheme \eqref{e:immpscheme} or \eqref{e:immpscheme2}  (computed with or without the Fixman correcting forces $\Dq V_{\M{fix},\nu}$), is of the order $2$ and consistent with respect to the original exact deterministic dynamics \eqref{e:langevin} (i.e., with $\gamma=0$, $\sigma=0$).
\epro
\begin{proof}
Following the proof of Proposition~\ref{p:shemesmallnu} we have, by the implicit function theorem, locally, the RATTLE scheme %
is a standard leapfrog scheme (see \cite{HaiLub02,LeiRei05}), and the associated local mapping depends continuously on $\nu^2$ when $\nu \to 0$. 
Therefore the usual calculation of the order of the leapfrog scheme (see \cite{HaiLub02})
 holds uniformly with respect to $\nu$, and \eqref{e:immpscheme} or 
\eqref{e:immpscheme2} are of order $2$ consistent uniformly in $\nu$. Now, the IMMP Hamiltonian  \eqref{e:immpscheme2} is a $\nu^2$ perturbation of the exact 
Hamiltonian, and the result then follows from applying a simple Gronwall argument.
\end{proof}

Similarly, we also obtain the limit towards the constrained/RATTLE scheme. 
\bpro[Large penalty limit] \label{p:shemelargenu}
Assume that the integrators \eqref{e:immpscheme} or \eqref{e:immpscheme2} are  locally well-defined, and that $\nu \to +\infty$ with $\frac{z_0}{\nu} \to \bar{z}_0$. Then the IMMP numerical scheme \eqref{e:immpscheme} or \eqref{e:immpscheme2} (computed with or without the Fixman correcting forces $\Dq V_{\M{fix},\nu}$) converges at order $\frac{1}{\nu^2}$ towards the constrained/RATTLE scheme for the rigid constraints $\xi(q_{n+1})  = \bar{z}_0$ (also computed with or without the associated geometric correcting Fixman forces $\Dq V_{\M{fix}}$ defined by \eqref{e:penV}).
\epro
\begin{proof}Following similar steps as in the proof of Proposition~\ref{p:shemesmallnu}, with auxiliary variables in the scaling $(\bar{z}, \bar{p}^z)=(\frac{z}{\nu}, \frac{p^z}{\nu})$,
and the following shift of Lagrange multipliers:
\[
\lambda \to \lambda +\frac{\dt}{2} \Dtwo U,
\]
 we obtain convergence at the order $\nu^{-2}$ to a decoupled scheme 
where the $(q,p)$ variables evolve according to the usual constrained/RATTLE scheme. Now from \eqref{e:penV}, $V_{\M{fix},\nu}$ is a perturbation 
of $V_{\M{fix}}$ of order $\nu^{-2}$ which completes the proof.
\end{proof}

\brem {\rm
We conclude this section with a practical recipe for tuning the mass matrix penalty. 
Identifying a suitable value of $\nu$ can be done, for instance, by computing the time averaged energy error \eqref{e:dyncritdt}  
in a dynamical simulation, or the Metropolis acceptance ratio \eqref{e:samplcritdt} in a sampling simulation, which gives a precise 
quantification of the time-step error.
Then increasing the penalty $\nu$ can save computational time as long as it leads to a reduction of the time-step error.
Indeed, this means that the selected fDOFs are limiting the time-step stability region. Prescribing the time-step error, 
a maximal time-step $\dt_{max}$ associated with the largest penalty $\nu_{max}$ that is able to improve stability can be obtained in this way.
Finally, one can set, for example, $\nu = \tfrac{\nu_{max}}{\dt_{max}} \dt$ in Scheme~\ref{e:immpscheme} to obtain an order two 
convergent scheme with an increased stability region.
}
\erem
\section{Numerical simulations for the $N$-alkane model}\label{s:numerexperiment}
The IMMP method is numerically tested on the united atom $N$-alkane model. 
The united atom model is a coarse-grained description of linear alkane isomers in which the hydrogen atoms
are not resolved and the molecule is modeled by a chain of particles interacting with effective potentials, see, e.g., \cite{MarSie98}. 
The integrator studied in the present section is described in Section~\ref{s:IMMP} and Section~\ref{s:numerinteg},
and is systematically compared with the exact dynamics, which is numerically 
integrated by a simple Leapfrog/Verlet
scheme with or without a thermostat, and with the constrained dynamics numerically integrated by the RATTLE scheme with or without a thermostat. 
Constraints are resolved using a simple Newton algorithm with a Gaussian linear solver. 
The Fixman forces are not resolved in the dynamical part, and thus the dynamics interpolates 
for a large penalty with
the ``rigidly'' constrained  dynamics. However, they are used in the Metropolis rule when sampling is considered.
 
The Metropolis/HMC step for the dynamics integrator as described by Scheme \ref{d:schemehmc} is added only when sampling is studied.

\begin{figure}[ht]
  \centerline{
  \includegraphics[width=13cm]{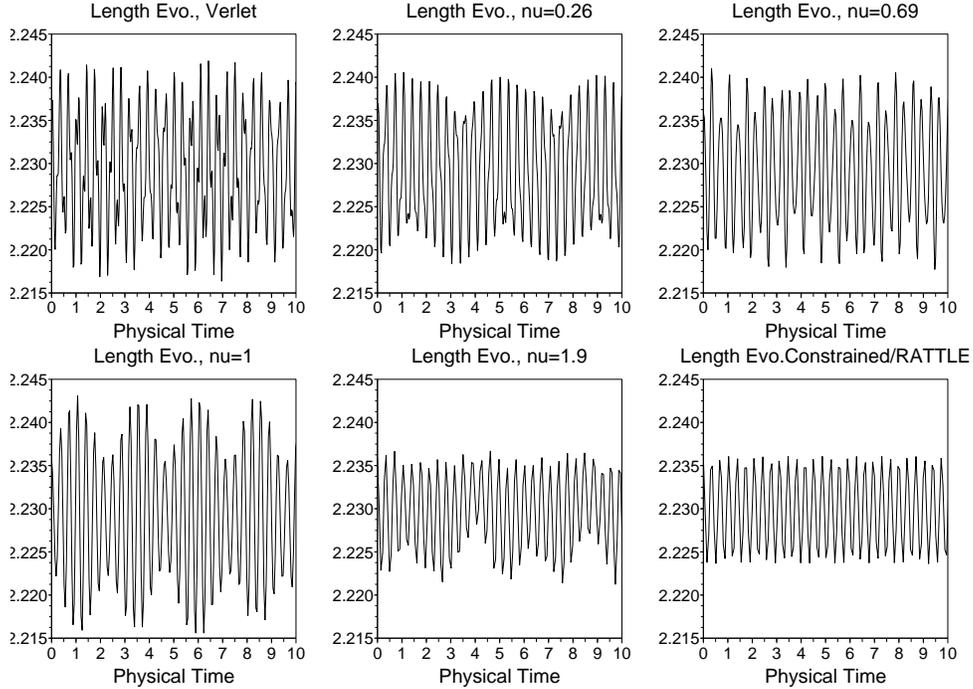} 
  }
\caption{\label{f:dyn_but_traj} Oscillations of the butane end-to-end length, for the Verlet integrator, 
the IMMP integrator, and the constrained integrator. 
Simulations are computed with a prescribed initial energy. Note the interpolation property of the IMMP scheme. 
This figure is associated to the frequency analysis in Figure~\ref{f:dyn_but_freq_dist}.}
\end{figure}

\begin{table}[h]
\begin{center}
\begin{tabular}{|c|cccccc|}
\hline
&Verlet  & $\nu = 0.5$ & $\nu = 1.0$ &  $\nu = 1.3$  & $\nu = 1.9$ &  RATTLE\\
\hline
$\dt^{\rm dyn}_c$ &  $.024 \, (1) $& $.032 \, (1) $& $.046 \, (1)$ &  $.059 \, (1)$ &  $.077 \, (1)$ & $.093 \, (2)$   \\
\hline
$\dt^{\rm sampl}_c$ & $ .013 \, (1) $ & $ .014 \, (1) $& $.022 \, (1)$ &  $.028 \, (1)$ &  $.035 \, (1)$ & $.049 \, (1)$   \\
\hline
\end{tabular}
   \caption{\label{t:dyn_but_tstep} Critical time steps $t_c^{\DYN}$ and $\dt^{\SAMPL}_c$ of the butane deterministic dynamics and 
    of the sampling schemes, respectively. Note the interpolation property from constrained dynamics to exact dynamics.}
\end{center}
\end{table}

\begin{figure}[ht]
  \centerline{
  \includegraphics[width=12cm]{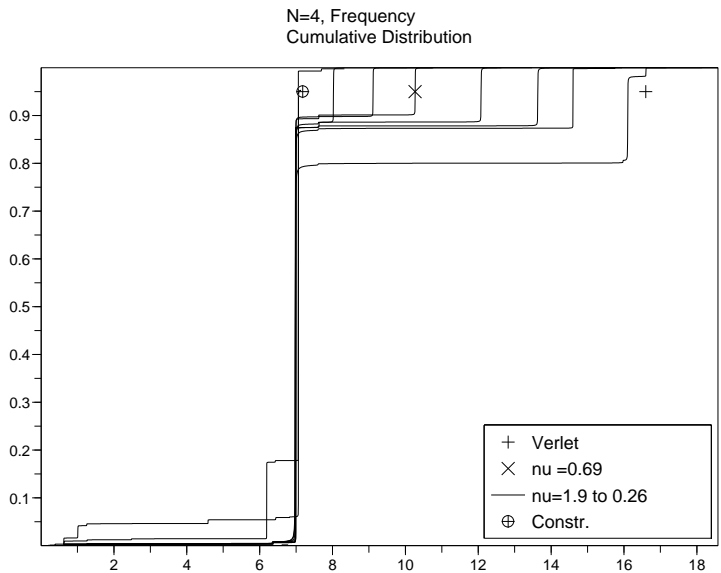}
}
\caption{\label{f:dyn_but_freq_dist} Frequency distribution of the end-to-end length oscillations of the butane Verlet dynamics and of the IMMP dynamics, 
cf. Figure~\ref{f:dyn_but_traj}. 
Note the two main frequencies (slow torsion and fast bond angles oscillations), and the slow components due to fluttering of resonances. }
\end{figure}

\begin{figure}[ht]
\centerline{
\includegraphics[width=12cm]{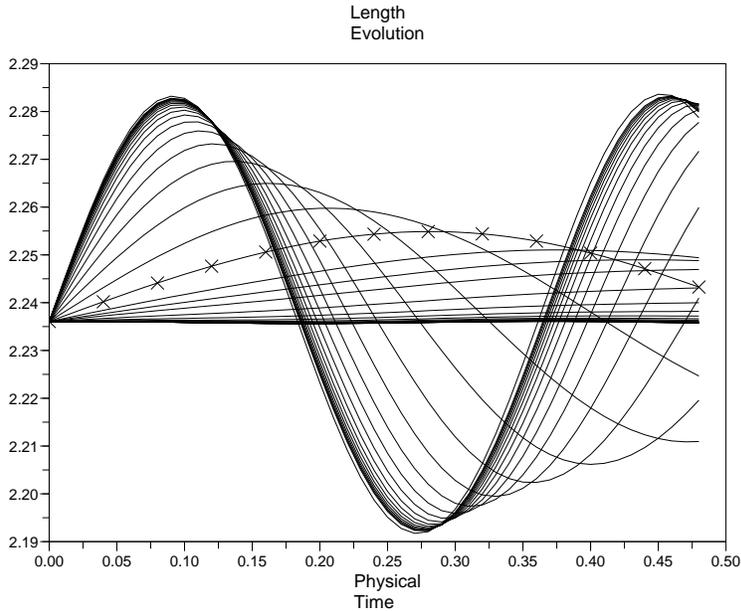}
}
\caption{\label{f:dyn_but_short} Interpolation of short trajectories of the butane length with IMMP dynamics, 
from constrained to Verlet dynamics. The marked trajectory corresponds to $\nu = 5.5$.
The Verlet dynamics is the oscillating limit.}
\end{figure}

\begin{figure}[ht]
\centerline{
\includegraphics[width=12cm]{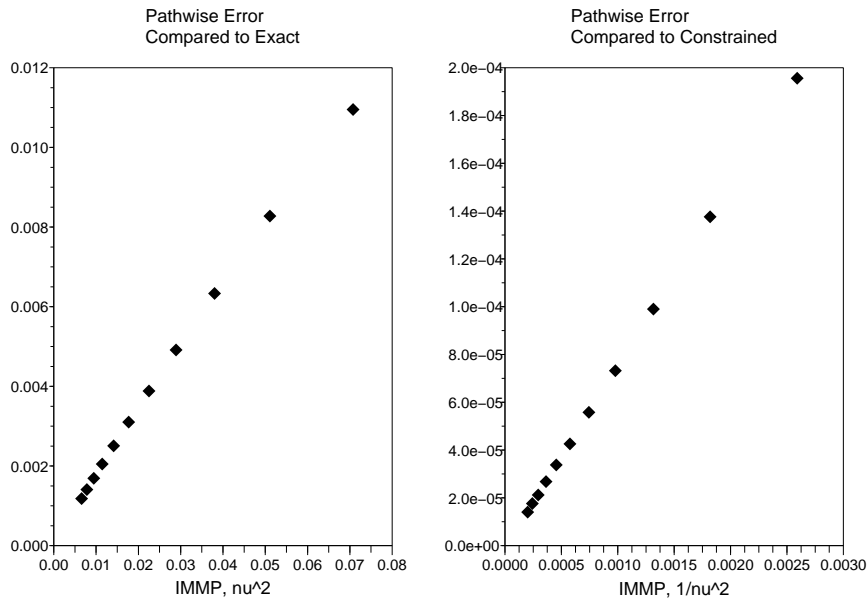}
}
\caption{\label{f:dyn_but_err} Convergence in terms of the $\ell^2$ pathwise error of the short time butane 
length IMMP dynamics compared to Verlet dynamics and constrained dynamics, with orders of convergence $\BIGO(\nu^2)$ and $\BIGO(\nu^{-2})$.}
\end{figure}

\begin{figure}[ht]
\centerline{
\includegraphics[width=12cm]{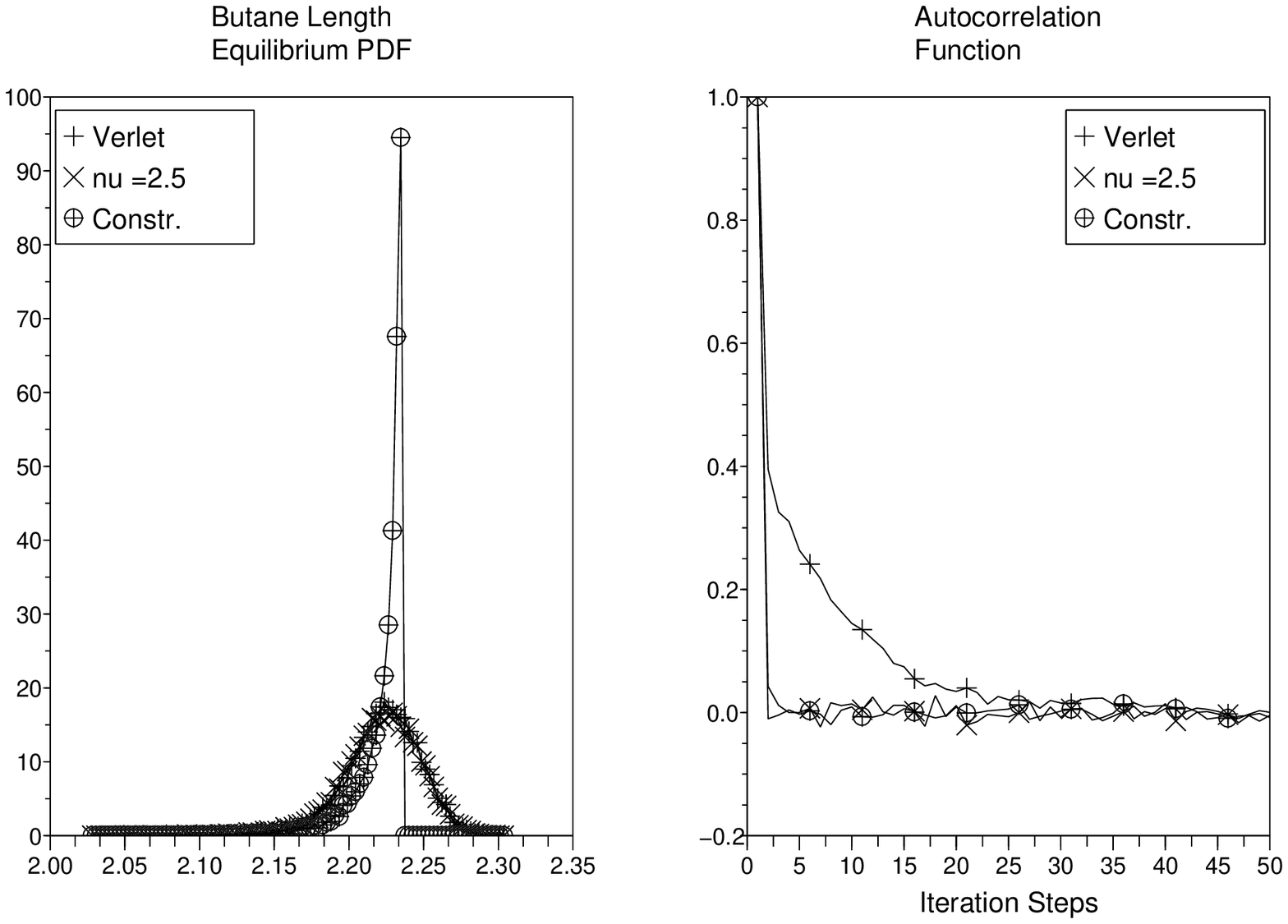}
}
\caption{\label{f:sampl_but} {\it Left:\/} 
Equilibrium PDF of the end-to-end length of the butane molecule with the GHMC scheme, using Verlet, IMMP
(penalty $\nu$), and constrained/RATTLE integrators. Note that the constrained integrator does not sample the correct measure. 
{\it Right:\/} The autocorrelation function in terms of iteration steps for the GHMC scheme comparing the IMMP and the Verlet integrator. 
The decrease in the $\ell^2$-decorrelation time \VIZ{decorr_time} is by the factor $1.8$.}
\end{figure}

\subsection{The $N$-alkane model}
The model consists of a chain of $N$ atoms with position vectors $q^{\PID{i}}\in\R^3$, hence the vector of DOFs is
$q=(q^{\PID{1}},\dots,q^{\PID{N}})\in\R^{3N}$. The mass of particles is normalized to be $m^\PID{i} = 1$ for all $i=1,\dots,N$.
The interaction potential $V$ consists of three short-range potentials that involve 2-body, 3-body, and 4-body terms
$$
 V(q) = \sum_{i,j\in\mathcal{I}_\BOND} V_\BOND(q^\PID{i},q^\PID{j}) + 
        \sum_{i,j,k\in\mathcal{I}_\ANGLE} V_\ANGLE(q^\PID{i},q^\PID{j},q^\PID{k}) +
        \sum_{i,j,k,l\in\mathcal{I}_\TORSION} V_\TORSION(q^\PID{i},q^\PID{j},q^\PID{k},q^\PID{l})\PERIOD 
$$

In the presented simulations we focus on short-range interactions as they are primarily responsible for the
stiffness of the system. Thus we omit the long-range Coulomb or Lenard-Jones
interaction potentials.

\smallskip\noindent{\it 2-body interactions.} The pair-wise bond potential $V_\BOND$ depends on the distance $r^\PID{i} = |q^\PID{i+1} - q^\PID{i}|$, 
$i=1,\dots,N-1$ of the bonded particles
$(q^\PID{i}, q^\PID{i+1})$ in the linear chain. Typically this potential is assumed to be harmonic with the equilibrium distance
$r^\PID{i} = r_0$. In all simulations considered here the bonds are treated as rigid bonds of the constant length $r_0=1$. Thus
defining the constrained DOFs
\[
 \xi_\BOND(q) = (r^\PID{1},\dots,r^\PID{N-1}) = (1,\dots,1)\PERIOD
\]

\smallskip\noindent{\it 3-body interactions.} The interaction of three consecutive particles  $(q^\PID{i},q^\PID{i+1},q^\PID{i+2})$ in the chain
is defined by the potential $V_\ANGLE(\theta^\PID{i})$, $i=1,\dots,N-2$  that depends on the angle between the vectors
$(q^\PID{i+2}-q^\PID{i+1})$ and $(q^\PID{i+1}-q^\PID{i})$. In the simulated model we assume the bending angle potential
\[
 V_\ANGLE(\theta) = \frac{A_0}{2}\sin^2(\theta - \frac{\pi}{2})\PERIOD
\]
Depending on simulation the bending angles can be mass-penalized or directly constrained leading to the definition
of penalized DOFs
\begin{eqnarray*}
 \xi_\ANGLE(q) = (\theta^\PID{1},\dots,\theta^\PID{N-2}) =& \frac{z_a}{\nu}\;\;\;&\mbox{in the case of mass-penalization,} \\
 \xi_\ANGLE(q) = (\theta^\PID{1},\dots,\theta^\PID{N-2}) =& 0 \;\;\;& \mbox{in the case of rigid constraints,}
\end{eqnarray*}
where $z_a\in\R^{N-2}$ are the auxiliary variables of the IMMP method.

\smallskip\noindent{\it 4-body interactions.} The last contribution $V_\TORSION(\phi^\PID{i})$
to the interaction potential depends  on the torsion (dihedral)
angles $\phi^\PID{i}$, $i=1,\dots,N-3$ which are defined as angles between two planes $\mathrm{span}\{(q^\PID{i+3}-q^\PID{i+2}),(q^\PID{i+2}-q^\PID{i+1})\}$ and 
$\mathrm{span}\{(q^\PID{i+2}-q^\PID{i+1}),(q^\PID{i+1}-q^\PID{i})\}$.  
We use a simple choice of the torsion potential
\[
 V_\TORSION(\phi) = - B_0 \cos(\phi)\PERIOD
\]
The torsion DOFs are mass-penalized with the auxiliary variable $z_b\in\R^{N-3}$
\[
 \xi_\TORSION(q) = (\phi^\PID{1},\dots,\phi^\PID{N-3}) = \frac{z_{b}}{\nu}\PERIOD
\]

\noindent{\it System of units.} For the purpose of computational tests we have chosen the system of units in which the
mass of united atoms is $m^\PID{1}=\dots =m^\PID{N}= m_0 = 1$, the equilibrium bond lengths are $r^\PID{1} = \dots = r^\PID{N-1} = r_0 = 1$
and the inverse temperature $\beta = 1/kT = 1$ at the ambient temperature $T = 300\,K$. 
The time step $\delta t = 0.01$ in these units corresponds to the physical time step $\delta \tilde t = \sqrt{\beta m_0} r_0 = 3\,\mathrm{fs}$.
Following physically relevant parameters (see, for example, \cite{Rap95,MarSie98}) with an artificially slightly stiffer bending angle potential 
leads to the angle potential constant $ A_0 = 500$ and the torsion potential constant $B_0 =20$.

\medskip

Applying direct constraints to bonds or bending angles in the presented formulation of the IMMP method
consists
in replacing the auxiliary variables $z$ associated with fully constrained DOFs by $0$. 
As described in Proposition~\ref{p:largepen} on the large penalty limit 
(see also Appendix~\ref{s:measures} and Section~\ref{s:highoscill} for the limit of stiff distributions), 
the associated dynamics is the standard ``rigid'' constrained dynamics. Then for sampling, the Fixman geometric corrector \VIZ{e:Vfixnu} 
is computed and accounts for the geometric bias of all the constrained DOFs including bonds. 
Note that the physical relevance of the geometric bias due to bonds is ambiguous since the quantum resolution may need to be introduced for
the bond interactions. Such a change can be easily incorporated by computing the Fixman corrector potential associated with the bonds only and subtracting 
it from the Fixman corrector potential related to all the constrained DOFs.

\subsection{Efficiency criteria}
The efficiency of the IMMP method as compared to the Verlet/Leapfrog integrator is quantified from two different viewpoints.

\smallskip
\noindent{\it Dynamics.} To compare computational efficiency of numerical schemes, a notion of critical time step has to be introduced. 
The critical time step $\dt^\DYN_c$ is defined implicitly through the formula
\begin{equation}\label{e:dyncritdt}
  \int \frac{\beta}{3N-N_c} [H_{n+1} - H_{n}]^+ \pen{\mu}(\pen{p}^n,q^n) = \ERR\COMMA
\end{equation}
where $N_c$ denotes the number of rigidly constrained DOFs, $[\cdot]^+$ is the positive part, $H_n$ is the energy of the numerical scheme 
at the $n$-th step, and $\pen{\mu}$ is the Boltzmann distribution associated with the Hamiltonian at hand as defined in \eqref{e:penboltz}. 
The quantity $\ERR$ is a prescribed typical non-dimensional error of the energy per degrees of freedom as compared to the inverse
temperature $\beta$. At least when $\ERR$ is small, this defines uniquely $\dt^{\DYN}_c$. Thus for a given $\ERR$, the larger $\dt^{\DYN}_c$ is, 
the less costly the method is in terms of force evaluations per integration time step.

The critical time step $\dt^{\DYN}_c$ is compared for the IMMP method and the Verlet/Leapfrog scheme. To achieve a fair comparison, it has then to be 
checked that the IMMP penalty, which introduces a dynamical modification, does not modify the relevant slow frequencies and slow varying components 
of the system on large time intervals. This can be done by comparing the distortion of the frequency distribution of a long trajectory. 
In the numerical tests we analyze for $t\in [0,T]$
the trajectory of the end-to-end molecule length $t\mapsto L(t)$, and we define the frequency density 
as the normalized square modulus of its Fourier transform $d(\omega) = |\hat{L}(\omega)|^2/Z_f$
where $Z_f$ is the normalization making $d(\omega)$ a probability density function, and $\hat{L}(\omega) = \int_0^T L(t) \exp(i \omega t)\, dt$.
We plot the cumulative distribution function
$$
D(\omega) = \int_0^\omega d(\omega')\,d\omega'\PERIOD
$$
Hence, in the figures dominant frequencies correspond to jumps in the associated cumulative distribution.
Of course, this comparison remains still largely qualitative.

\smallskip
\noindent{\it Sampling.} When using a Metropolis correction step in the numerical simulation, as in Scheme~\ref{d:schemehmc}, the critical 
time step is correlated to the rejection rate of the Metropolis step, since the larger the former is, the more rejections will occur. 
Thus for sampling methods, the time step $\dt^{\SAMPL}_c$ is tuned in order to achieve a given rate of rejection $\RATE$ satisfying
\begin{equation}\label{e:samplcritdt}
 \int \exp\pare{-\beta [H_{n+1} - H_{n}]^+} \pen{\mu}(\pen{p}^n,q^n) = 1-\RATE,
\end{equation}
where $\exp\pare{-\beta [H_{n+1} - H_{n}]^+}$ is the Metropolis weight that appears in Generalized Hybrid Monte-Carlo methods, as defined 
in Scheme~\ref{d:schemehmc}. 
Then, computational efficiency is defined by comparing $\dt^{\SAMPL}_c$ with the mixing time \emph{in terms of physical time} of the resulting Markov Chain. 
Equivalently, and more directly, the results are presented as the mixing time \emph{in terms of iteration steps}. Such comparison requires choosing
an appropriate notion of the mixing time.
At least when the momenta are overdamped (i.e., $\gamma \to +\infty$), the chain becomes reversible, and the mixing time can be 
rigorously defined using the spectral gap of the underlying Markov kernel. In the present work, the decorrelation time of relevant observables 
is used. Assuming the initial state $q_0$ of the system is at equilibrium, the normalized autocorrelation function associated with 
a given position observable $q \mapsto \phi(q)$ is 
\[
C_n(\phi) := \frac{\E[\phi(q_{n}) \phi(q_0)]}{\E[\phi(q_{0})^2]}\COMMA
\]
and it can be computed in large time simulations using path averages. 
Then $\ell^2$-decorrelation time $n_{\mathrm{corr}}$, calculated in terms of iteration steps, is given by the formula
\begin{equation}\label{decorr_time}
 n_{\mathrm{corr}} = 2\sum_{n=0}^{+\infty} C_n(\phi)^2\PERIOD
\end{equation}
This quantity corresponds to the approximate number of steps of the chain for the observable $\phi$ to decorrelate 
from its past values.

\subsection{Numerical results}
For comparisons we choose an often studied observable defined as the end-to-end distance of
the alkane chain.

\medskip
\noindent{\sc I. The butane model ($N=4$).} The bonds between atoms are rigidly constrained, and bond angles are treated using the 
IMMP method with Scheme~\ref{e:immpscheme}.

     \smallskip\noindent{\it Dynamical behavior.} The dynamics of the butane model is studied using 
     deterministic dynamics with a prescribed initial energy. 
     As shown in Figure~\ref{f:dyn_but_traj}, the IMMP dynamics with different penalty $\nu$ yield  
     an interpolation from the mixed 
     torsion/bond angles oscillations of the exact dynamics, to the simple torsion oscillation of 
     the constrained dynamics. 
     Depending on the frequency introduced by the mass penalization, some fluttering resonance can 
     be observed in Figure~\ref{f:dyn_but_traj}  between the torsion and the bond angles.

     Behavior observed in Figure~\ref{f:dyn_but_traj} is related to the frequency density of the 
     end-to-end length oscillations which are depicted in Figure~\ref{f:dyn_but_freq_dist}.
     The bond angles oscillation frequency appears to slow down with the IMMP penalization 
     and is eventually setting
     to a single frequency of the constrained dynamics. 
     Note the slow modes introduced by the resonances. 

     The dynamical interpolation of the IMMP from the exact dynamics to the constrained dynamics is 
     demonstrated on short time trajectories in Figure~\ref{f:dyn_but_short}. The associated 
     convergence orders, $\BIGO(\nu^2)$ and $\BIGO(\nu^{-2})$,
     respectively, are captured in Figure~\ref{f:dyn_but_err}. 
     The time step stability is studied in Table~\ref{t:dyn_but_tstep}. The IMMP method enables an  
     increase of the critical time step 
     of the Verlet scheme, following the interpolation property. 

    \smallskip\noindent{\it Sampling behavior.} Exact sampling of the equilibrium distribution on a 
    very large time scale, 
    whatever the value of the IMMP penalization, 
    is shown in Figure~\ref{f:sampl_but}. The distribution of the butane length for constrained 
    bond angles is clearly distorted. 
    The mixing time to equilibrium is also studied. The autocorrelation function of the length 
    evolution in terms of iteration steps is plotted
    in Figure~\ref{f:sampl_but}, and the faster convergence of the IMMP method is demonstrated. 
    The $\ell^2$-decorrelation time \VIZ{decorr_time} is decreased by the 
    factor $1.8$ using the IMMP method. 

\begin{figure}[ht]
\centerline{
\includegraphics[width=12.0cm]{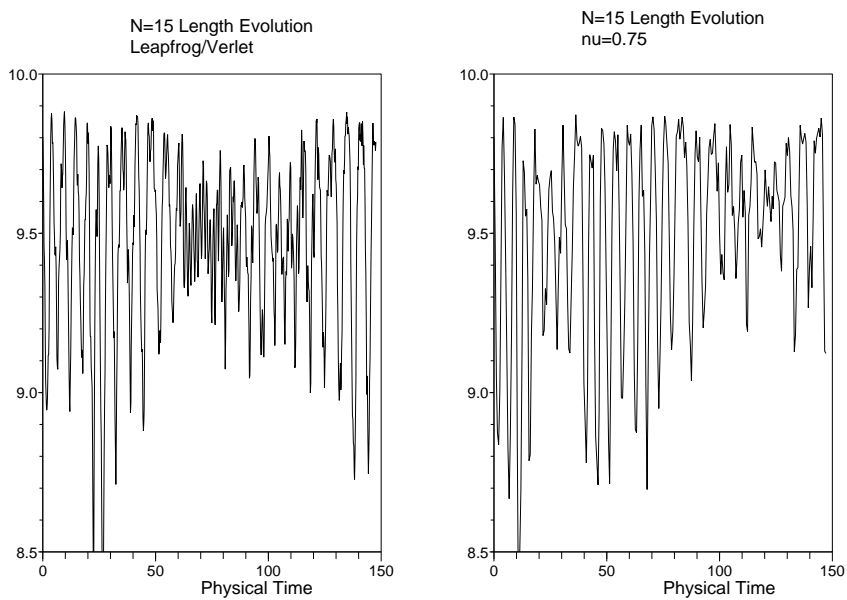}
}
\caption{\label{f:dyn_alk_traj} The trajectory of the $N=15$-alkane dynamics for the Verlet scheme and the IMMP scheme. 
Note that the IMMP penalty does not modify substantially the low frequencies/slowly varying components. The frequency analysis is presented in Figure~\ref{f:dyn_alk_freq_dist}}
\end{figure}

\begin{figure}[ht]
\centerline{
\includegraphics[width=12.0cm]{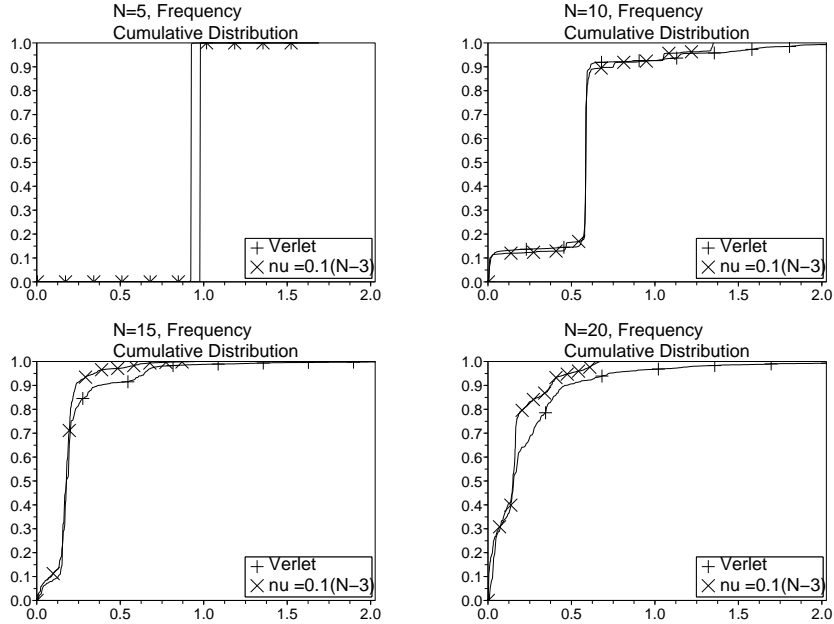}
}
\caption{\label{f:dyn_alk_freq_dist} Spectral densities of the end-to-end length of the alkane for $N=5,10,15,20$. 
Note that the IMMP penalty does not modify substantially the low frequencies/slowly varying components.}
\end{figure}

\begin{figure}[ht]
\centerline{
\includegraphics[width=12.0cm]{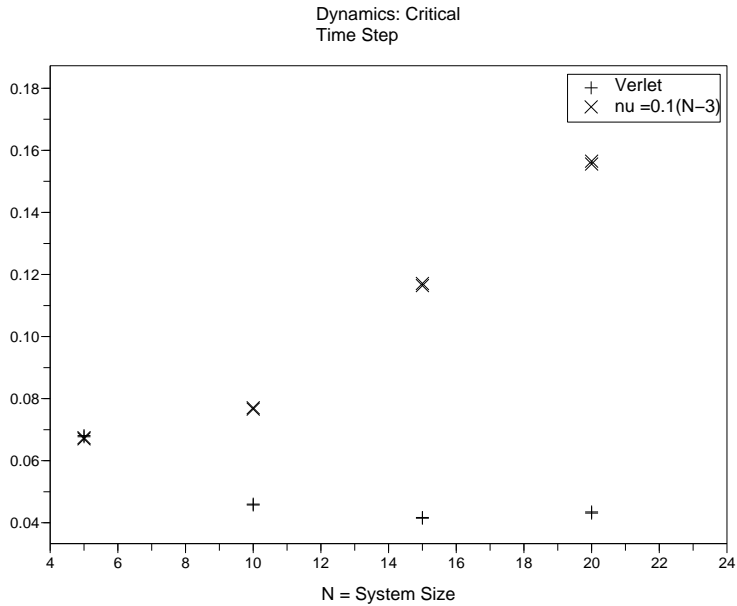}
}
\caption{\label{f:dyn_alk_tstep} Critical time steps, with error bars, $\dt_{c}^{\DYN}$ of the Verlet and IMMP dynamics with respect to the system size.}
\end{figure}

\begin{figure}[ht]
\centerline{
\includegraphics[width=12.0cm]{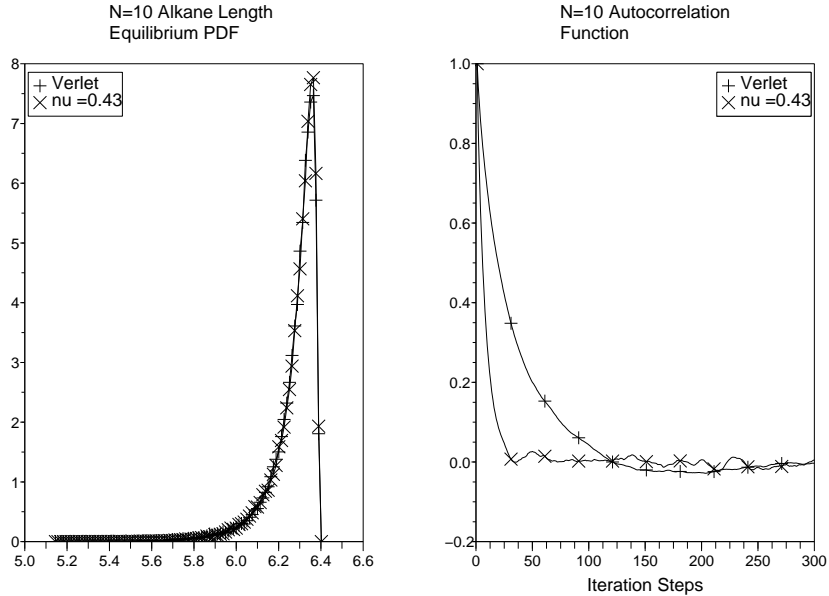}
}
\caption{\label{f:sampl_alk_dist} {\it Left:} Equilibrium PDF of the end-to-end length of the alkane chain with the GHMC scheme, 
  using Verlet and IMMP integrator (penalty $\nu$). 
  {\it Right:} The autocorrelation function in terms of iteration steps for the GHMC scheme comparing the IMMP and the Verlet integrator. 
  The system size is $N=10$}
\end{figure}

\begin{figure}[ht]
\centerline{
\includegraphics[width=12.0cm]{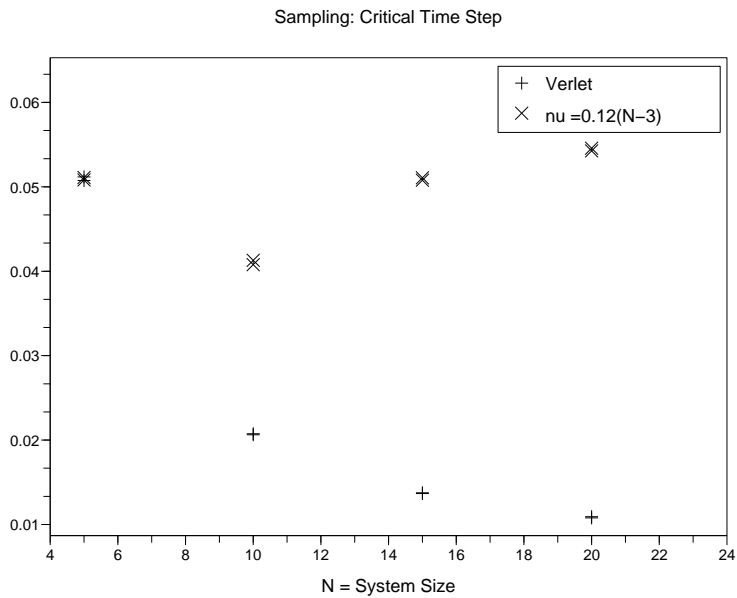}
}
\caption{\label{f:sampl_alk_tstep} Critical time steps, with error bars, $\dt_{c}^{\SAMPL}$ of the Verlet and IMMP dynamics with respect to the system size. 
    Note that the IMMP penalty heals the degeneracy of the rejection rate for large systems.}
\end{figure}

\begin{figure}[ht]
\centerline{
\includegraphics[width=12.0cm]{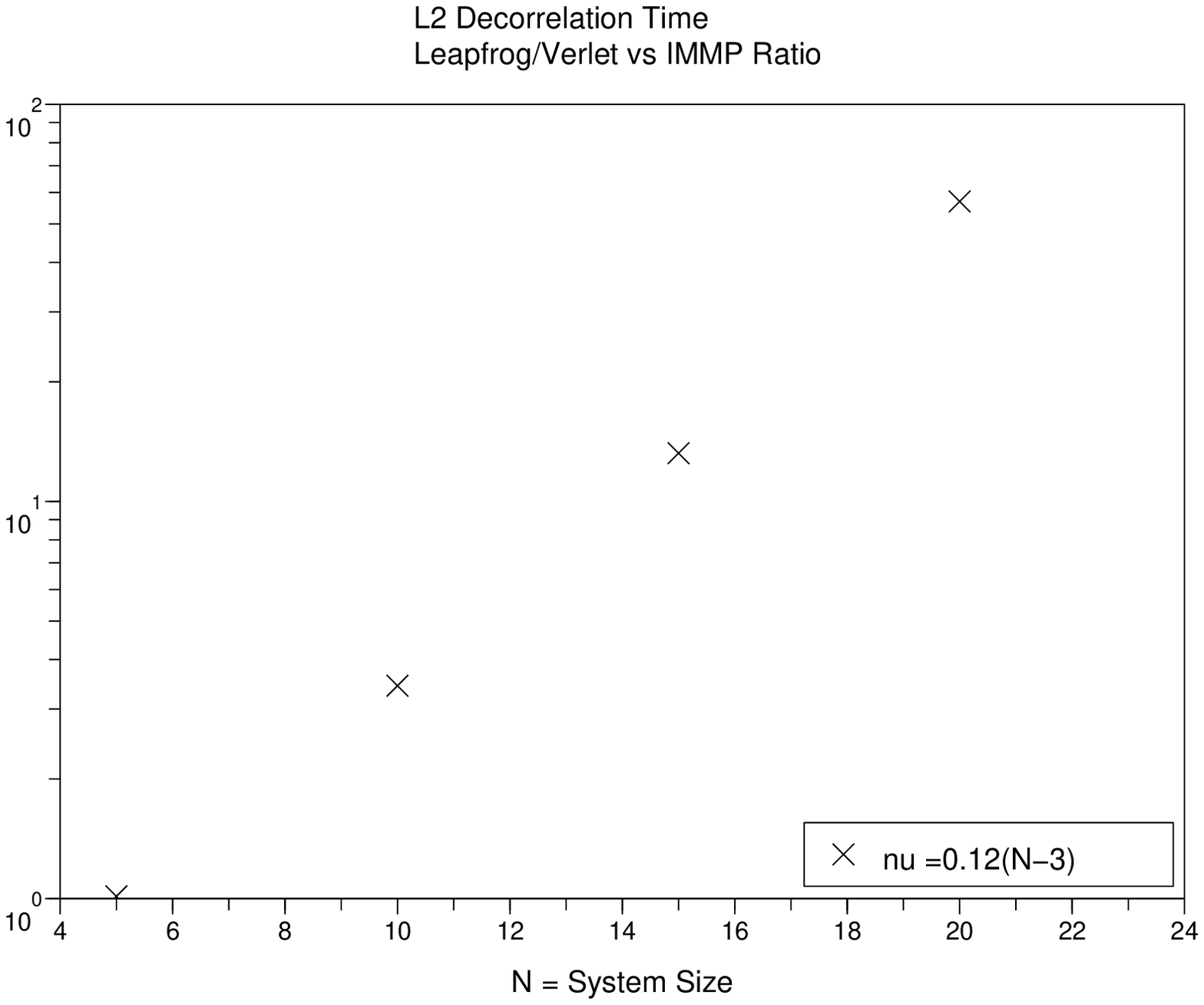}
}
\caption{\label{f:sampl_alk_timeratio} Comparison the $\ell^2$-decorrelation time of the end-to-end alkane length in terms of Monte-Carlo iteration steps. 
  The ratio between the case of the Verlet integration and the case of the IMMP integration is depicted. Note that the vertical axis is in the logarithmic scale 
  and an exponential gain occurs. Error bars are plotted, but do not appear at this scale.}
\end{figure}

\medskip
\noindent{\sc II. The alkane model ($N=5,\dots,20$).}
In this test the bonds and bending angles between atoms are rigidly constrained, and torsion angles are treated using the IMMP method. 
Note that the rigid constraints on torsions would lead to a rigid molecule, losing completely the evolution of the molecule length. 
The IMMP penalty is increased with the system size using the linear scaling $\nu = \bar {\nu}N$.
In Section~\ref{s:Nanal} we present a mathematical justification of this scaling for linear systems. 
In principle, introducing inertia in the torsion angles adds weights on the diagonal of the mass-tensor in the internal coordinates of the
molecule. This can reduce the lowest eigenvalues of the mass-tensor, and can explain the reduction of the multi-scale nature of the system 
that are observed in the simulations below. 

        \smallskip\noindent{\it Dynamical behavior.} The frequency of oscillations of the alkane length 
         in Figure~\ref{f:dyn_alk_traj}, and in the spectral analysis 
         in Figure~\ref{f:dyn_alk_freq_dist} are not substantially modified by the IMMP penalization. One can observe a small group of fast oscillations 
         in the middle of the Verlet dynamics plot in Figure~\ref{f:dyn_alk_traj} which is not present in the IMMP case. 
         This translates in the top of the spectral plots in Figure~\ref{f:dyn_alk_freq_dist} where a cut-off of the fastest oscillatory scales 
         for the IMMP case occurs. 

         The critical time steps $\dt^{\rm dyn}_c$ with respect to the system size are depicted in Figure~\ref{f:dyn_alk_tstep}. 
         The gain in time stepping increases more than linearly with the system size $N$.
         This behavior, however, depends on the initial scaling $\nu = \bar {\nu}N$ %
         of the IMMP penalty.

        \smallskip\noindent{\it Sampling behavior.} Exact sampling of the equilibrium distribution on very large time scales
          (whatever the value of the IMMP penalization) is shown for $N=10$ in Figure~\ref{f:sampl_alk_dist}, with the auto-correlation of 
          the latter observable with the gain in mixing time of the IMMP dynamics.  
          The precise ratio of the $\ell^2$-decorrelation time between the IMMP integrator and the Verlet one is given in Figure~\ref{f:sampl_alk_timeratio} 
          for different system sizes. It increases again more than linearly in $N$, in fact exponentially for this particular system. 
          The associated critical time steps $\dt^{\SAMPL}_c$ are in Figure~\ref{f:sampl_alk_tstep}.
          We observe that the critical time step increases with large $N$ which demonstrates that in the present case, the 
          IMMP method heals the decrease of the Metropolis rejection rate of for large systems (see also \cite{IzaHam04}).

\medskip
\noindent{\it Conclusions.}
The presented numerical studies demonstrate that for integrating the dynamics, the IMMP allows for relaxation of the time-step stability restrictions.
In the case of sampling methods the IMMP method decreases the decorrelation time measured in terms of Monte-Carlo iteration steps leading to
more efficient sampling algorithm.
In both cases the improvements are increasing with system size $N$.

\section{The infinite stiffness limit}\label{s:highoscill}
Throughout this section, one introduces a potential function with an explicit dependence with respect to the fast variables $(q,z) \mapsto U(q,z)$ together with a stiffness parameter
$\epsilon$. The potential energy $V$ can then be written in the form
\begin{equation*}
V(q ) = U(q,\frac{\xi(q)}{\epsilon}) \COMMA\;\;\;\mbox{with a confining assumption}\;\;
\inf_{q\in \R^d} U(q,z) \geq K(z) \COMMA
\end{equation*}
where $K(z) > c\log|z|$ as $|z|\to\infty$.
The functions $\xi$ are then indeed ``fast'' degrees of freedom (fDOFs) in the limit $\epsilon \to 0$, 
the system being confined to the slow sub-manifold $\SMAN = \set{q \SEP \xi(q)=0}$.

We prove, that under appropriate scaling of the mass penalty $\nu_{\epsilon} = \tfrac{\nusc}{\epsilon}$, 
the IMMP method is asymptotically stable in the stiff limit, converging towards standard effective dynamics on the slow manifold $\SMAN$.

\subsection{Thermostatted stiff systems}
The canonical distribution becomes
\be \label{e:boltzmannstiff}
   \mu_{\epsilon}(dp\,dq) = \frac{1}{Z_{\epsilon}}\EXP{ - \beta (\frac{1}{2} p^{T} \MASST^{-1} p +
   U(q,\frac{\xi(q)}{\epsilon})) } dp\, dq \PERIOD
\ee
In the infinite stiffness limit ($\epsilon \to 0$) the measure concentrates on the slow manifold $\SMAN$.
The limit is computed using the co-area formula (see Appendix~\ref{s:measures} for relevant definitions of surface measures).
In order to characterize the limiting measure we introduce the effective potential
\be \label{e:Veff}
V_{\M{eff}}(q) = -\frac{1}{\beta} \ln \int \EXP{- \beta U(q,z) } \, dz \PERIOD
\ee
\blem\label{l:boltzlim}%
In the infinite stiffness limit ($\epsilon \to 0$),
the highly oscillatory canonical distribution \eqref{e:boltzmannstiff} converges $\mu_\epsilon \WEAKLY \mu_0$
(in distribution) towards $\mu_{0}(dp\,dq)$,
 which is supported on $\SMAN$, and defined as
\be \label{e:boltzlim}
\mu_{0}(dp\,dq) = \frac{1}{Z_{0}} \EXP{-\beta (\frac{1}{2} p^{T} \MASST^{-1} p +V_{\M{eff}}(q))} \,dp \, \delta_{\xi(q) = 0}(dq) \PERIOD
\ee
Its marginal distribution in position is given, up to the normalization, by
\be \label{e:boltzlimmarg}
 \EXP{- \beta V_{\M{eff}}(q)}\delta_{\xi(q) = 0}(dq) \PERIOD
 \ee
\elem
\begin{proof}
It is sufficient to consider distributions in the position variable $q$ only. Let $\mathcal{U}^\delta$ be a
$\delta$-neighborhood of $\SMAN$ where
$dq =\epsilon^{\nc} \delta_{\xi(q) = \epsilon z}(dq) \, dz$.
We construct a decomposition $\ph = \ph_1 + \ph_2$ of continuous bounded observables such that $\supp \ph_1\subset\mathcal{U}^\delta$
and $\supp \ph_2 \cap \mathcal{U}^{\delta/2} = \emptyset$.
Using the confining property of $U(q,\cdot)$ we obtain
\[
\int \ph(q) \EXP{-\beta U(q,\frac{\xi(q)}{\epsilon})} dq =\epsilon^{\nc}
\int \ph_1(q) \EXP{-\beta U(q,z)} \delta_{\xi(q)= \epsilon z}(dq) \, dz  + \BIGO(\EXP{-\beta K(\delta/2\epsilon)})\PERIOD
\]
By continuity of $\epsilon \mapsto \int \ph_1(q) \EXP{-\beta U(q,z)} \delta_{\xi(q)= \epsilon z}(dq)$ and by
the dominated convergence theorem
\[
  \int \ph_1(q) \EXP{-\beta U(q,z)} \delta_{\xi(q)= \epsilon z}(dq) \, dz \to \int \ph_1(q)
  \EXP{- \beta V_{\M{eff}}(q)}\delta_{\xi(q) = 0}(dq)=\int \ph(q) \EXP{- \beta V_{\M{eff}}(q)}\delta_{\xi(q) = 0}(dq)  \COMMA
\]
and the result follows after normalization.
\end{proof}

The infinite stiffness limit ($\epsilon \to 0$) of highly oscillatory dynamics has been studied in a series of papers
\cite{RubUng57, Tak80, Kam85, BorSch97, Rei95, Rei00}. The limiting dynamics can be fully characterized in special cases.
For example, when the highly oscillatory potential is linear and non-resonant
(at least almost everywhere on the trajectory, see \cite{Tak80}), it can be described through adiabatic effective potentials. 
See also \cite{JahLub06, LebLeg07} for a recent work on some related numerical issues.
However, when the system is thermostatted, one can postulate an ``ad hoc'' effective dynamics (\cite{Rei00})
exhibiting the appropriate limiting canonical distribution given by \eqref{e:boltzlimmarg}.
Such dynamics can be obtained by constraining the system to
the slow manifold $\SMAN$, and adding  a correcting entropic potential, sometimes called Fixman corrector from \cite{Fix78}, 
which is due to the geometry of $\SMAN$, and is  given by
\be\label{e:Vfix}
V_{\M{fix}}(q) = \frac{1}{2\beta} \ln \pare{ \mathop{\det} G(q) } \COMMA
\ee
where $G(q)$ is the $\nc \times \nc$ Gram matrix defined in \eqref{e:Gram}.

In general, since the effective potential \eqref{e:Veff} is not explicit, one may need to couple the system with virtual fast degrees
of freedom to enforce the appropriate effective dynamics associated with \eqref{e:Veff}.
The resulting extended Hamiltonian is then defined on the state space $T^{*} \left( \SMAN \times \R^{\nc} \right)$
(the cotangent bundle) and is given by

\begin{equation} \label{e:Heff}
\syst{
&H_{ \M{eff} }(p,p_{z},q,z) =  \frac{1}{2} p^{T} \MASST^{-1} p  +\frac{1}{2} p_{z}^{T} \MASSTZ^{-1} p_{z}
+ U(q,z) + V_{\M{fix}}(q) & \\
&\xi(q) = 0 \PERIOD & \TAG{C}
}
\end{equation}

\begin{Def}[Effective Langevin process with constraints]\label{d:consteff}
The constrained Langevin process associated
with Hamiltonian \eqref{e:Heff} is defined by the following stochastic differential equations
  \begin{equation}\label{e:consteff}
  \syst{
    &\dot{q} =  \MASST^{-1} p &   \\
    &\dot{z} =  \MASSTZ^{-1} p_{z} &  \\
    &\dot{p} = -\Done U(q,z) - \Dq V_{\M{fix}}(q) -\gamma \dot{q} + \sigma \dot{W} - \Dtwo \xi\, 
                \dot{\lambda} &   \\
    &{\dot{p}_{z}} = -\Dtwo U (q,z)-\gamma_{z} \dot{z} + \sigma_{z} \dot{W}_{z} &\\
    &\xi(q)  = 0\COMMA                 & \TAG{C}
  }\end{equation}
where $\Done$ and $\Dtwo$ are respectively derivatives with respect to the first and second variable of the function $U(q,z)$, $\dot{W}$ (resp. $\dot{W}_{z}$ ) is the standard multi-dimensional white noise, $\gamma$ (resp. $\gamma_{z}$)
  a $d\times d$ (resp. $\nc\times \nc$) symmetric positive semi-definite dissipation matrix,
  $\sigma$ (resp. $\sigma_{z}$) is the fluctuation matrix satisfying $\sigma \sigma^{T} =\tfrac{2}{\beta} \gamma$
  (resp. $\sigma_{z} \sigma_{z}^{T} =\tfrac{2}{\beta} \gamma_{z}$).
  The processes $\lambda \in \mathbb{R}^{\nc}$ are Lagrange multipliers associated with the constraints
  $(C)$ and adapted with respect to the white noise.
\end{Def}

We formulate reversibility of this process as a separate lemma.
\blem
The process defined in \VIZ{e:consteff} is reversible with respect to the associated canonical distribution whose marginal distribution in $(q,p)$ variables is
\begin{equation}\label{e:mueff}
\mu_{\M{eff}}(dp \, dq) = \frac{1}{Z_{\M{eff}}}\EXP{-\beta \pare{ \frac{1}{2} p^{T} \MASST^{-1} p + V_{\M{eff}}(q) + V_{\M{fix}} (q)} } \sigma_{T^{*}\SMAN}(dp\,dq)
\end{equation}
with the $q$-marginal
\[
 \EXP{- \beta V_{\M{eff}}(q)}\,\delta_{\xi(q) = 0}(dq) \PERIOD
\]
When $\gamma$ and $\gamma_{z}$ are strictly positive definite, the process is ergodic.
\elem
\begin{proof}
The process \VIZ{e:consteff} is a Langevin process with mechanical constraints, exhibiting reversibility properties
with respect to the associated Boltzmann canonical measure (see the summary in Appendix~\ref{s:langevinapp}).
Then the $q$-marginal is obtained by remarking that the integration of any function of
$\tfrac{1}{2} p^{T} \MASST^{-1} p + \tfrac{1}{2} p_{z}^{T} \MASSTZ^{-1} p_{z}$ with respect to $d p_z \, \sigma_{T^{*}_q\SMAN}(dp)$
results in a constant independent of $q$.
\end{proof}

The properties of thermostatted highly oscillatory systems are summarized in Table~\ref{t:sum}.
\begin{table}[H]
\begin{center}
\begin{tabular}{|c|cccc|}
\hline
&   \multicolumn{2}{c|}{Finite stiffness} &      \multicolumn{1}{c|}{ Infinite stiffness limit }    &   \multicolumn{1}{c|}{ Infinite stiffness  }            \\
&   \multicolumn{2}{c|}{ $\epsilon >0$  }  &      \multicolumn{1}{c|}{$\epsilon \to 0$}  & \multicolumn{1}{c|}{$\epsilon = 0$}                       \\
\hline
Dynamics & Highly & & \multicolumn{1}{c|}{Adiabatic }      & Effective with \\
        &  oscillatory  &                             & \multicolumn{1}{c|}{(if non-resonant)}    &     constraints \\
         &  $+$ fluct./diss.   &                               &\multicolumn{1}{c|}{ $+$ non-Markov fluct./diss.}   &   $+$ fluct./diss.            \\

\hline
Statistics & Canonical  &  &  \multicolumn{1}{c|}{Positions on $\mathcal{M}_{0}$,} & Canonical on $T^{*}\mathcal{M}_{0}$,  \\
           &            &                                &  \multicolumn{1}{c|}{free velocities.}     &  geometric corrector. \\
\hline
Numerics & Leapfrog/Verlet  &      &    \multicolumn{1}{c|}{Time-step  }                   &      Leapfrog/Verlet with     \\
     &     $+$ fluct./diss.       &                 &    \multicolumn{1}{c|}{restrictions ($\dt = o(\epsilon)$)}               &       constraints     \\
      &         &                      &    \multicolumn{1}{c|}{}               &        $+$ fluct./diss.           \\
\hline

\end{tabular}
\caption{\label{t:sum} Stiff  Hamiltonian systems and associated commonly used numerical methods
($\SMAN$ denotes the slow manifold).
Two different schemes are required
for the stiff system and its effective Markovian approximation.}
\end{center}
\end{table}

\subsection{Stability of the IMMP dynamics}\label{s:infstiff}
We assume that the mass-matrix penalty parameter $\nu\equiv\nueps$ grows to infinity in such a way that
$\lim_{\epsilon\to 0} \epsilon \nu_\epsilon = \nusc$.

The original Hamiltonian with the stiffness parameter is expressed explicitly as
\be\label{e:Heps}
H_{\epsilon}(p,q) =
\frac{1}{2} p^{T} \MASST^{-1} p    + U(q,\frac{\xi(q)}{\epsilon}) \COMMA
\ee
and including the mass-matrix penalization one gets
\be\label{e:Hnueps}
	\peneps{H}(\peneps{p},q) = \frac{1}{2} \peneps{p}^{T} \peneps{\MASST}^{-1} \peneps{p}
                   + U(q,\frac{\xi(q)}{\epsilon}) + V_{\M{fix},\nueps}(q)  \COMMA
\ee
or in its implicit formulation
\begin{equation}\label{e:Himmpeps}\syst{
   &H_{\M{IMMP}}(q,z,p,p_{z}) = \frac{1}{2} p^{T} \MASST^{-1} p  + \frac{1}{2} p_{z}^{T} \MASSTZ^{-1} p_{z} +
   U(q,\frac{z}{\nueps\epsilon} ) + V_{\M{fix},\nueps}(q)\COMMA & \\
   &\xi(q) = \frac{1}{\nueps} z\PERIOD & \TAG{C_{\nueps}}
}\end{equation}
One immediately sees that $H_{\M{IMMP}}$ is non-singular when $\epsilon\to 0$ and converges to the effective Hamiltonian on the slow manifold,
\begin{equation}\label{e:Himmpstab}
\syst{
     &  H_{\M{eff},\nusc}(q,z,p,p_{z}) = \frac{1}{2} p^{T} \MASST^{-1} p  + \frac{1}{2} p_{z}^{T} \MASSTZ^{-1} p_{z}
                                          + U(q,\frac{z}{\nusc}) + V_{\M{fix}}(q) &\\
       &\xi(q) = 0 \PERIOD &\TAG{C} \PERIOD
}
\end{equation}
The expression \eqref{e:Himmpeps} represents a minor generalization of  $H_{\M{eff}}$ in \eqref{e:Heff},
but it leads to the same canonical marginal distribution $\mu_{\M{eff}} (dp \, dq)$ in $(p,q)$ variables as given by \eqref{e:mueff}.
The continuity in $\epsilon$ of $H_{\M{IMMP}}$ implies stability of the associated dynamics and their numerical integrators.
We  first derive the limits of the original and penalized canonical distribution.
\bpro[Limits of canonical distributions]
Consider the canonical distributions $\peneps{\mu}(dp \, dq$ associated with the mass penalized Hamiltonian \eqref{e:Hnueps}, but considered with respect to the 
variables $(p\equiv\MASST\peneps{\MASST}^{-1} \peneps{p},q)$. 
In the sense of weak convergence of measures we have 
$\peneps{\mu}\WEAKLY \mu_{\mathrm{eff}}$ as $\epsilon\to 0$ with
$\mu_{\M{eff}}$ defined by \eqref{e:mueff}.
\epro
\begin{proof}
The first convergence is  proved in  Lemma~\ref{l:boltzlim}. For the second one,  the following notation will be used
\[
\delta_{q,\epsilon z}(dq)  = \delta_{\xi(q) = \epsilon z}(dq)\COMMA\;\;\;\mbox{and}\;\;
\delta_{p,\epsilon p_z}(dp)  =\delta_{p^T \MASST^{-1} \nabla \xi(q)  = \epsilon\MASSTZ^{-1} p_z}(dp) \PERIOD
\]

To prove the convergence towards $\mu_{\mathrm{eff}}$
we consider a $\delta$-neighborhood $\mathcal{U}^\delta$ of $\SMAN$ where
\begin{eqnarray*}
  dp \, dq &=&\frac{\epsilon^{2\nc}}{\DET\MASSTZ} \delta_{q,\epsilon z}(dq) \, dz \;
  \delta_{p,\epsilon p_z}(dp) \, dp_z \COMMA \\
\end{eqnarray*}
and a decomposition of the bounded observable (in $(p,q)$ variables) $\ph = \ph_1 + \ph_2$
such that $\supp \ph_1\subset\mathcal{U}^\delta$ and $\supp\ph_2\cap\mathcal{U}^{\delta/2} = \emptyset$.
Thus, keeping in mind that $\peneps{p} = \peneps{\MASST} \MASST^{-1} p$, and using the confining property of the potential $U(q,\cdot)$ we obtain
\begin{equation}\label{testfun}
\int \ph(p,q) \EXP{-\beta \peneps{H}} d \peneps{p} \, dq = \int \ph_1(p,q) \EXP{-\beta \peneps{H}} d \peneps{p} \, dq +
 \BIGO(\EXP{-\beta K(\delta/\epsilon)}) \equiv I_\epsilon +  \BIGO(\EXP{-\beta K(\delta/2\epsilon)}) \PERIOD
\end{equation}
Applying the change of variables $\peneps{p} = \peneps{\MASST} \MASST^{-1} p $ yields
\[
d\peneps{p} = \DET(\peneps{\MASST }\MASST^{-1})\, dp
= \nueps^{2\nc} \,\DET{\MASSTZ} \,\DET(G+\frac{1}{\nueps^2}\MASSTZ^{-1})\, dp \COMMA
\]
and setting $  \epsilon\MASSTZ^{-1} p_z = \nabla_q \xi\,\MASST^{-1} p $ and
$\epsilon z = \xi(q) $ we get
\[
\peneps{H}(\peneps{p},q) = \frac{1}{2} p^{T} \MASST^{-1} p    + \nueps^2 \epsilon^2
p_z^{T} \MASSTZ^{-1} p_z + U(q,\frac{z}{\nueps\epsilon}) + V_{\M{fix},\nueps}(q)
= H_{\M{IMMP}}(q,z,p,p_z)\PERIOD
\]
Thus substituting back to \VIZ{testfun} we obtain
\[
 I_\epsilon = (\nueps \epsilon)^{2 \nc} \int \ph_1 \EXP{-\beta H_{\M{IMMP}}(q,z,p,p_z)} \DET(G+\frac{1}{\nueps^2}\MASSTZ^{-1}) \, \delta_{p,
\epsilon p_z}(dp) \, dp_z\, \delta_{q,\epsilon z}(dq) \, dz \COMMA
\]
and thus
\[
I_\epsilon \xrightarrow[\epsilon \to 0]{}  \nusc^{2\nc}
\int \ph_1 \EXP{-\beta H_{\M{eff},\nusc}(q,z,p,p_z)}  \DET(G) \, \delta_{p,\epsilon p_z}(dp) \, dp_z \,
\delta_{\xi(q) = 0}(dq) \, dz
\PERIOD
\]
Using the co-area formula we obtain
\[
\DET(G)\, \delta_{\nabla \xi(q)\MASST^{-1} p  = 0}(dp)\delta_{\xi(q) = 0}(dq) = \sigma_{T^* \SMAN}(dp \, dq) \COMMA
\]
which leads to the final result after integration of the $(p_z,z)$ variables and normalization.
\end{proof}
\brem {\rm
Due to the fast oscillations, the distribution of impulses in the limiting distribution $\mu_{0}$ in
\eqref{e:boltzlim} is uncorrelated, whereas after the mass-matrix penalization, the limiting distribution \eqref{e:penboltz}
has almost surely co-tangent impulses (i.e., satisfying the constraints $\nabla_q \xi \MASST^{-1} p =0$).
This explains the role of the corrected potential energy $V_{\M{fix}}$ taking into account the curvature of $\mathcal{M}_{0}$.
}
\erem

In the next step we inspect the infinite stiffness asymptotic of the penalized dynamics.
\bpro[Infinite stiffness limit]
	When $\epsilon \to 0$ with $\nu\equiv\nueps \sim \tfrac{\nusc}{\epsilon}$ and
        $V(q,\xi(q))= U(q,\tfrac{\xi(q)}{\epsilon})$,
        the IMMP Langevin stochastic process \eqref{e:IMMP} converges weakly towards
        the following coupled limiting processes with constraints
	\begin{equation}\label{e:limdyn}
        \syst{
    	  &\dot{q} =  \MASST^{-1} p \COMMA &  \\
    	  &\dot{p} = -\Done U(q,\frac{z}{\nusc}) - \Dq V_{\M{fix}}(q) -\gamma \dot{q} + \sigma \dot{W} - \Dq\xi \dot{\lambda}\COMMA  & \\
    	  &\xi(q)  =  0  \COMMA                                        &  \TAGG{C} \\
    	  &\dot{z} =  \MASSTZ^{-1} p_{z} \COMMA &  \\
    	  &{\dot{p}_{z}} = -\frac{1}{\nusc}\Dtwo U(q,\frac{z}{\nusc})-\gamma_{z} \dot{z} + \sigma_{z} \dot{W}_{z} \PERIOD &
	}
        \end{equation}
	where $\Done$ and $\Dtwo$ are respectively derivatives with respect to the first and second variable of the function $U(q,z)$,
        and $\{\lambda_t\}_{t\geq 0}$ are adapted stochastic processes defining the Lagrange multipliers associated with the constraints $(C)$.

	The process $\{q_t,p_t\}_{t\geq 0}$ defines an effective dynamics with constraints
        (Definition~\ref{d:consteff}) for thermostatted highly oscillatory systems.
        It is reversible with respect to its stationary canonical distribution given by $\mu_{\mathrm{eff}}$ \eqref{e:mueff},
        and is ergodic when $(\gamma,\gamma_z)$ are strictly positive definite.
\epro
\begin{proof} The proof is similar to the proof of Proposition~\ref{p:largepen}. Here we have
\[
 \Dq U = \Done U + \frac{1}{\epsilon} \Dq \xi \Dtwo^T U \COMMA
\]
and \eqref{e:IMMP} translates, up to a change of Lagrange multipliers, into
\begin{equation}\syst{
    &\dot{q} =  \MASST^{-1} p  &  \\
    &\dot{z} =  \MASSTZ^{-1} p_{z} &  \\
    &\dot{p} = -\Done U - \Dq V_{\M{fix},\nu}(q) -\gamma \dot{q} + \sigma \dot{W} - \Dq \xi \,\dot{\lambda} & \\
    &{\dot{p}_{z}} = - \frac{1}{\nueps \epsilon} \Dtwo U-\gamma_{z} \dot{z} + \sigma_{z} \dot{W}_{z} +  \frac{\dot{\lambda}}{\nueps} & \\
    &\xi(q) = \frac{z}{\nueps}\COMMA & \TAG{\peneps{C}} \PERIOD
  }\end{equation}
The rest follows the proof of Proposition~\ref{p:largepen}.

\end{proof}
\brem {\rm
When $\nusc \to +\infty$, by a classical averaging argument (see, e.g., \cite{Kam85}), one can check that the limiting dynamics
are the effective dynamics pointed out in \cite{Rei00}
\begin{equation} \label{e:avlimdyn}\syst{
    	&\dot{q} =  \MASST^{-1} p & \\
    	&\dot{p} = -\Dq U_{\M{eff} } (q) - \Dq V_{\M{fix}}  -\gamma \dot{q} + \sigma \dot{W} - \Dq\xi\, \dot{\lambda}& \\
    	&\xi(q)  = 0  \PERIOD & \TAG{C}
}\end{equation}
with the stationary canonical distribution \VIZ{e:mueff}.
}
\erem

\subsection{Stability of the IMMP integrator}
The numerical scheme (Scheme~\ref{d:scheme} proposed for the IMMP method \eqref{e:IMMP}) is also stable in the limit of infinite stiffness
$\epsilon \to 0$.
Recall that we consider a reversible, measure preserving numerical flow
$\Phi_{\dt}^{\nueps}(p,p_z,q,z)$ associated with Hamiltonian \eqref{e:Himmpeps}  $H_{\M{IMMP}}$ with constraints (modified potentials could similarly be considered).
\bpro[Asymptotic stability]
In the limit $\epsilon\nueps \to \nusc$, the numerical flow $\Phi_{\dt}^{\nueps}$ associated with the leapfrog/Verlet integrator with constraints
for the IMMP Hamiltonian \eqref{e:Himmpeps} converges towards the numerical flow
$\Phi_{\dt}^{\nusc}$, which is the leapfrog/Verlet integrator with geometric constraints associated with effective Hamiltonian \eqref{e:Himmpstab} on the slow manifold.
\epro
\begin{proof}
The statement is a direct consequence of the implicit function theorem and  the continuity of the leapfrog integrator
with constraints \eqref{e:immpscheme} with respect to the parameter $\nusc=\epsilon \nueps$.
Indeed, considering the shift of Lagrange multipliers
$\lambda \to \lambda +\tfrac{1}{\epsilon} \Dtwo U$
and taking the limit $\epsilon \to 0$ we obtain the appropriate leapfrog scheme
\begin{equation*}\syst{
		&p_{n+1/2} = p_{n} -  \frac{\dt}{2} \Done U(q_{n},\frac{z_n}{\nusc}) - \Dq \xi(q_{n})  \lambda_{n+1/2} &\\
		&p_{n+1/2}^z = p_{n}^z -  \frac{\dt}{2 \nusc} \Dtwo U(q_{n},\frac{z_n}{\nusc})     & \\
		&q_{n+1} =  q_{n}+ \dt \MASST^{-1} p_{n+1/2}            & \\
		&z_{n+1} = z_n + \dt \MASSTZ^{-1} p_{n+1/2}^z             & \\
		&\xi(q_{n+1}) = 0                       & \TAG{C_{1/2}} \\
		&p_{n+1} = p_{n+1/2} -  \frac{\dt}{2} \Done U(q_{n+1},\frac{z_{n+1}}{\nusc}) - \Dq \xi(q_{n+1})  \lambda_{n+1} &\\
		&p_{n+1}^z = p_{n+1}^z -  \frac{\dt}{2 \nusc} \Dtwo U(q_{n+1},\frac{z_{n+1}}{\nusc})  & \\
		&\Dq \xi(q_{n+1})\MASST^{-1}p_{n+1} = 0 \PERIOD & \TAG{C_1}
}\end{equation*}
\end{proof}

By convergence of the Hamiltonian \eqref{e:Himmpeps} to \eqref{e:Himmpstab},  similar asymptotic stability properties holds when a Metropolis step is introduced.

The results and properties discussed in this section are summarized in Table~\ref{t:sum2}.
\begin{table}[H]
\begin{center}
\begin{tabular}{|c|cccc|}
\hline
& Zero mass    &  \multicolumn{2}{|c|}{Positive }          &    \multicolumn{1}{|c|}{ Infinite }  \\
& penalization &  \multicolumn{2}{|c|}{ mass-penalization} &    \multicolumn{1}{|c|}{ stiffness limit }  \\

&$\nu = 0$              &  \multicolumn{2}{|c|}{ $\epsilon,\nu >0$  } &    \multicolumn{1}{|c|}{ $\epsilon\to 0$, $\tfrac{\nu}{\epsilon} \to \nusc$}               \\
\hline
Dynamics & Highly oscillatory & IMMP && Effective with  \\
         & $+$ fluct./diss.   &  $+$ fluct./diss.     &                                &    constraints$+$ fluct./diss.    \\

\hline
Statistics  &Canonical&   Canonical with      &                &  Canonical on \\
            &         & correlated velocities &                &      $T^* \SMAN$                 \\
\hline
Numerics & \multicolumn{4}{c|}{ IMMP $+$ fluct./diss. }     \\
\hline
\end{tabular}

\caption{\label{t:sum2} The IMMP dynamics and the Verlet numerical integration are both asymptotically stable in the infinite stiffness
regime if $ \tfrac{\nu}{\epsilon} \to \nusc < +\infty$.
If the mass-penalization vanishes ($\nu =0$) one recovers the original physical stiff system.
The canonical distribution is always exact in the position variable. Notice that due to the penalized mass-matrix ($\nu>0$)
the statistics have correlated velocities.}
\end{center}
\end{table}

\section{Numerical analysis of a harmonic particle chain}\label{s:Nanal}
In this section we present rigorous analysis for a special case of the linear chain with {\it harmonic} interactions. The analysis
supports scaling properties of the IMMP, with respect to the size of the chain, observed in numerical simulations of the
general linear alkane chains.
We consider the thermodynamic limit $N\to +\infty$ where $N$ is the size of the system. It is shown that the macroscopic dynamics of the IMMP method behaves continuously (uniformly with $N$, and in the $L_2$ norm for the position profile) with respect to the re-scaled mass penalty parameter $\nusc$.

At the same time, the time-step stability of the IMMP numerical scheme \eqref{e:Hnum} is compared with the
standard Verlet scheme, and the critical time step is shown to be increased by a factor $\nusc N$.

From the spectral point of view, the IMMP method behaves in this linear case as a low-pass filter. This proves, in this simplified case, the ability of IMMP method to respect macroscopic dynamical equivalence, while saving computational time up to a factor of order $\BIGO(N)$.

\subsection{Conservation of macroscopic dynamics}

The model we consider consists of a chain of particles which interact through the harmonic (quadratic) potential $v_{\INTP}(r) = r^2/2$.
Each particle is also individually submitted to a macroscopic confining exterior (quadratic) potential
$v_{\EXTP}(r)$.
After converting to the non-dimensional form the typical quantities involved in the model enable us to write
a scaling at the mass-transport level where the dynamics of the chain is described by the Hamiltonian
\be\label{e:HN}
H_{N}(q,p) = \frac{1}{2} p^{T}  p + \sum_{i=1}^{N-1}  v_{\INTP}(\dg_{i} q ) + \sum_{i=1}^{N} v_{\EXTP}(q_{i})  \COMMA
\ee
and by a coupling with an exterior thermal bath at the re-scaled inverse temperature $\beta_N = \beta N^{-1}$.
In the expression \eqref{e:HN} the functions $r \in\R \mapsto v_{\INTP}(r)\in\R$ and $q \mapsto v_{\EXTP}(q)\in\R$
are the smooth interaction potential
and the exterior potential, respectively.
The linear operator $\dg:\mathbb{R}^{N} \to \R^{N-1}$, having the components
 \[
 \dg_{i} q = \frac{q_{i+1}-q_{i}}{1/N}\COMMA\;\;\;i=1,\dots,N-1\COMMA
 \]
represents the discrete gradient associated to the chain with the Neumann boundary conditions.
Its transpose operator is denoted  $(\dg)^T:\mathbb{R}^{N-1} \to \R^{N}$. The discrete Laplace operator
is then defined as $\dl = -(\dg)^T\dg$.
The particles are represented by their re-scaled positions $q=(q_1,...,q_N)$, so that the typical position
and deviation of  $q$ is formally of order $1$ with respect to $N$.
This can be seen by considering particles in the chain as indexed by $x= \tfrac{i}{N} \in [0,1]$.
We choose to work with such scaling in $N$ that it prescribes the macroscopic
timescale of the chain profile at order one with respect to $N$.
Following our general construction we obtain the mass-penalized
Hamiltonian
\be\label{e:HNpen}
\penN{H}(\penN{p},q) = \frac{1}{2} \penN{p}^{T} \MASST^{-1}_{\nu_N} \penN{p}
                       + \sum_{i=1}^{N}  v_{\INTP}( \dg_i q )
                       + \sum_{i=1}^{N} v_{\EXTP}(q_{i}) \PERIOD
\ee
We chose the penalizing matrix to be the identity matrix $\MASSTZ = \ID$, hence the penalized
mass-tensor becomes $\penN{\MASST} = \ID - \nusc^{2}\dl$,
and the fluctuation/dissipation tensor is taken proportional to the identity matrix.
The system of stochastically perturbed equations of motions
then becomes
\be \label{e:spdelin}
\syst{
& \dot{q} =   ( \ID - \nusc^{2} \dl  )^{-1} p & \\
&\dot{p} =   \dl\,q  - v'_{\EXTP}(q) - \gamma \dot{q} +   \sigma \sqrt{N}\dot{W}  \COMMA&
}
\ee
with fluctuation/dissipation identity $\sigma^2 = 2\beta^{-1}\gamma$.
The associated canonical equilibrium distribution is then given by the re-scaled inverse temperature $\beta_N = \beta N^{-1}$.

In order to treat the limit $N\to\infty$, we introduce the $\ell_2$-norm in the position space
\[
\norm{q}_{\ell_2}^2 := \frac{1}{N} \sum_{i=1}^N q_i^2  = \frac{1}{N} q^T q \COMMA
\]
as well as the $h_{-1}$-norm in the momentum space
 \[
 \norm{p}_{h_{-1}}^2 = \norm{(-\dl)^{-1/2}(p-\frac{1}{N}\sum_{i=1}^N p_i)}^2_{\ell_2}+
    \left(\frac{1}{N}\sum_{i=1}^N p_i\right)^2.
 \]
In the above expression, $\tfrac{1}{N}\sum_{i=1}^N p_i$ can be seen as the orthogonal projection in $\ell^2$ on 
the one dimensional kernel of the Neumann discrete Laplacian $\dl$.  The quadratic form $\norm{q}_0^2 +\norm{p}_{-1}^2$ 
endows the phase-space with a Hilbert space structure.
\bpro[Convergence of the macroscopic dynamics]
Assume that the exterior potential $v_{\EXTP}$ is bounded and that its derivative satisfies the Lipschitz condition
$$
\norm{v'_{\EXTP}(q_2) - v'_{\EXTP}(q_1)}_{h_{-1}} \leq L_v \norm{q_2 - q_1}_{\ell_2} \COMMA
$$
where $L_v$ is independent of $N$. For any $T>0$, let $t \mapsto (p^{\nusc}(t),q^{\nusc}(t))$ be the solution, for $t\in [0,T]$, of the evolution equation \VIZ{e:spdelin} with the initial condition
\[
(p^{\nusc}(0),q^{\nusc}(0)) = (\penN{\MASST}^{-1/2} p^0(0),q^0(0)) \COMMA
\]
where $(p^0(0),q^0(0))$ is distributed according to the original equilibrium canonical distribution (associated 
with \eqref{e:HN} and $\beta_N = \beta N^{-1}$). Then for all $t\in [0,T]$ one has the uniform convergence
 \[
 \lim_{\nusc \to 0} \limsup_{N \to + \infty } \EXPECT{\norm{q^{\nusc}(t)-q^{\nusc=0}(t)}^2_{\ell_2}} = 0 \PERIOD
 \]
\epro
\begin{proof} We write $X = (q,p)$,  and introduce the norm
\[
\norm{X}_{\nusc} = \norm{q}_{\ell_2} + \norm{ \penN{\MASST}^{-1/2} p }_{h_{-1}}\PERIOD
\]
The system \eqref{e:spdelin} becomes a stochastic differential equation in the form
\begin{equation}\label{eqoper}
d X^{\nusc}_t = A_{\nusc} X^{\nusc}_t + F(X^{\nusc}_t) + \Sigma d W_t \COMMA
\end{equation}
where by definition
\[
A_{\nusc} = \bmat
               0                         & (\ID -\nusc^2 \dl )^{-1}  \\
               \dl & 0 \emat \COMMA\;\;
\;\;
F(X) = \bmat
       0 \\-v'_{\EXTP}(q)-  \gamma p
       \emat
        \COMMA\;\;
\Sigma =  \bmat
       0 \\
       \sqrt{N}\sigma
       \emat\PERIOD
\]
Duhamel formula gives an implicit expression for differences of solutions of \VIZ{eqoper} with the same noise
\begin{eqnarray}\label{e:spdemild}
X^{\nusc}_{t} - X^{0}_{t} &=& \pare{\EXP{A_{\nusc} t} - \EXP{A_{0} t}} X^{0}_0 + \int_{0}^{t}\pare{\EXP{A_{\nusc} (t-s)} - 
\EXP{A_{0} (t-s)} }\,(F(X^{0}_s)\,ds + \Sigma \,dW_s)  \nonumber \\
&& + \EXP{t A_{\nusc} } (X^{\nusc}_0 - X_{0}^0) +
                  \int_{0}^{t}\EXP{A_{\nusc} (t-s)}( F(X^{\nusc}_{s})-F(X^0_{s})) ds \PERIOD
\end{eqnarray}
We estimate the individual terms on the right hand side  in \VIZ{e:spdemild}.
We define $P$ as the coordinate transformation associated with the orthonormal spectral decomposition
\[
-\dl = P^{-1} \DIAG(\delta_{0},\dots,\delta_{N-1}) P\COMMA
\]
where $P P^{T} = \ID$. The eigenvalues of the discrete Neumann Laplacian are given, for $k=0,\dots,N-1$, by
\begin{equation}\label{eigenval}
\delta_{k} =
4N^2\sin^2\left(\frac{k\pi}{2N}\right) \limop{\sim}_{N \to \infty} k^2 \pi^2 \PERIOD
\end{equation}
Denoting the spectral coordinates
 \[
\hat{X}=(\hat{q},\hat{p})=(N^{-1/2}Pp, N^{-1/2}Pq)
\]
we have
\[
 \norm{X}^2_{\nusc} = \hat{p}^2_0 + \sum_{k=1}^{N-1} \frac{\delta_k}{1+\nusc^2\delta_k} \hat{p}^2_k +\sum_{k=0}^{N-1} \hat{q}^2_k \PERIOD
\]
The spectral decomposition leads to a block diagonal form of the operator $\EXP{A_{\nusc} t}$ with diagonal
$2\times 2$ blocks in the spectral basis
\[
\EXP{\widehat{A_{\nusc}}^{(0)} t} = \bmat
        1 & t  \\
          0 &  1
         \emat \COMMA
\]
as well as for $k=1,\dots,N-1$
\[
\EXP{\widehat{A_{\nusc}}^{(k)}} = \bmat
         0  &  (1 + \nusc^{2}\delta_k)^{-1}\\
\delta_k         &  0
         \emat \COMMA
\]
where $\widehat{A_{\nusc}}^{(k)}$ is the $2\times 2$ block associated with the coordinates $(\hat{q}_k,\hat{p}_k)$. Since $\widehat{A_{\nusc}}^{(k)}$ conserves the $k$-mode energy $\delta_k  (\hat{q}_k)^2 + (1+\nusc^2 \delta_k)^{-1} (\hat{p}_k)^2$, one can check that for any $N \geq 1$ the operator norm
\[
\normop{\EXP{A_{\nusc} t}}_{\nusc}^2 \leq  2+2t^2 \PERIOD
\]
Similarly, since in the sense of symmetric matrices $\penN{\MASST}^{-1/2} \leq \ID$, we have the bound
\[
\norm{ F(X^{\nusc}) - F(X^0) }_{\nusc} \leq \norm{\penN{\MASST}^{-1/2} \pare{v'_{\EXTP}(q^{\nusc}) -v'_{\EXTP}(q^0) +\gamma (p^{\nusc}-p^0) }}_{h_{-1}} \leq (L_F + \gamma)\norm{X^{\nusc} - X^0 }_{\nusc} \PERIOD
\]
Using independence of Brownian increments we compute
\begin{eqnarray*}
 \EXPECT{\norm{   \int_{0}^{t} \pare {\EXP{A_{\nusc} (t-s)} -\EXP{A_{0} (t-s)}}\Sigma \,dW_s }^2_{\nusc}} &= &\int_{0}^{t}  \sum_{i=1}^{N} \norm{ (\EXP{A_{\nusc} (t-s)}-\EXP{A_{0} (t-s)}) \Sigma_{.,i}}^2_{\nusc} ds \PERIOD
\end{eqnarray*}
Applying Gronwall lemma in \VIZ{e:spdemild} and collecting all terms we obtain
\begin{equation}\label{e:gronwallN}
\EXPECT{\norm{ X^{\nusc}_{t} - X^{0}_{t}}_{\nusc} ^2} \leq C_T \pare{ \EXPECT{\norm{ X^{\nusc}_{0} - X^{0}_0}_{\nusc} ^2
               + m_T }} \COMMA
\end{equation}
where $C_T$ is independent of $N$, and with $X^0$ being distributed canonically $m_T$ is given
\begin{eqnarray*}
 m_T &=& \sup_{t\in[0,T]} \pare{ \EXPECT{\norm{ \pare{ \EXP{A_{\nusc} t} - \EXP{A_{0} t}} X^0 }_{\nusc}^2}
         + \EXPECT{\norm{ \pare{ \EXP{A_{\nusc} t} - \EXP{A_{0} t}} F(X^0) }_{\nusc}^2}
         + \sum_{i=1}^N  \norm{ \pare{ \EXP{A_{\nusc} t} - \EXP{A_{0} t}} \Sigma_{.,i} }^2_{\nusc} }\PERIOD
\end{eqnarray*}
For a given random vector $X$ such that $\EXPECT{\norm{X}_0^2} < +\infty$, Parseval identity
and the inequality $\norm{\cdot}_{\nusc} \leq \norm{\cdot}_0$ imply
\begin{equation}\label{e:specserie}
\EXPECT{\norm{ \pare{ \EXP{A_{\nusc} t} - \EXP{A_{0} t}} X }_{\nusc}^2} =
  \sum_{k=1}^{N-1} \EXPECT{ \norm{ \pare{\EXP{\widehat{A_{\nusc}} t} - \EXP{\widehat{A_{0}}t} } \hat{X} }_{k,\nusc}^2}
  \leq 2 \sum_{k=1}^{N-1} \EXPECT{\norm{  \hat{X} }_{k,0}^2} \COMMA
\end{equation}
where $\norm{\cdot }_{k,\nusc}$ is the restriction to the $k$-th mode $(\hat{q}_k,\hat{p}_k)$.
Then one has, by orthogonality of $P$,
\[
 \sum_{i=1}^N \EXPECT{\norm{  \widehat{\Sigma_{.,i}}  }_{k,0}^2} =
      \sigma^2  \sum_{i=1}^{N} P^2_{k,i} \frac{1}{\delta_k} \leq \frac{\sigma^2}{\delta_k} \PERIOD
\]
Up to normalization, the distribution of $X^0$ has the density 
$\EXP{-  \tfrac{\beta}{N} \sum_{i=1}^d v_{\EXTP}(q^0_i)  }$ with respect to the Gaussian 
distribution with the covariance matrix $\beta^{-1} \ID$ for momenta variables,
and the covariance matrix $(\beta \dl)^{-1}$ for positions. Thus we have the bound
\[
 \EXPECT{\norm{  \widehat{X^0} }_{k,0}^2} \leq 2\EXP{4\beta \norm{v_{\EXTP}}_{\infty}} \frac{1}{\delta_k \beta} \COMMA
\]
as well as
\[
 \lim_{N \to \infty } \EXPECT{\norm{\widehat{F(X^0)} }_{k,0}^2}
     \leq \EXP{4\beta \norm{v_{\EXTP}}_{\infty}}
     \EXPECT{\norm{ \mathcal{F} \circ F \circ \mathcal{F}^{-1}( \widehat{G^0})}_{k,0}^2} \COMMA
\]
where $\mathcal{F}$ denotes the Fourier series expansion on $[0,1]$ with Neumann conditions, and
$(\widehat{G^0}_k)_{k \geq 1}$ are canonical centered Gaussian i.i.d. variables with the covariance matrix
$\beta^{-1}\bmat 1 & 0\\ 0 & \tfrac{1}{k^2\pi^2} \emat$. By the Lipschitz assumption
the series is bounded
\[
\sum_{k=1}^{+\infty} \EXPECT{\norm{ \mathcal{F} \circ F \circ \mathcal{F}^{-1}( \hat{G^0})}_{k,0}^2}
    \leq (L_v + \gamma) \EXPECT{\norm{\hat{G^0}}_{0}^2} = \sum_{k=1}^{+\infty} \frac{2(L_v + \gamma)}{\beta k^2\pi^2}.
\]
Since $\lim_{\nu \to 0}\normop{\EXP{\widehat{A_{\nusc}}^{(k)} t} - \EXP{\widehat{A_{0}}^{(k)} t}} = 0$,
one can take the limit $N \to +\infty$ and use the uniform convergence of the series in \eqref{e:specserie} to obtain
$\lim_{\nusc \to 0} \lim_{N \to +\infty} m_T =0$ in \eqref{e:gronwallN}. The convergence of the initial condition $\lim_{\nusc \to 0} \lim_{N \to +\infty} \E \norm{ X^{\nusc}_{0} - X^{0}_0}_{\nusc}^2 $ follows by using\ similar arguments. The proof is complete.
\end{proof}

\subsection{Relaxation of time-step stability restriction}
To demonstrate improved stability properties of time integration algorithms we consider the IMMP scheme \eqref{e:immpscheme} associated with the mass-matrix
penalized Hamiltonian \eqref{e:HNpen}. Note that when the constraints are linear, the leapfrog scheme (RATTLE) applied to an implicit
Hamiltonian is identical to the usual leapfrog scheme for the associated explicit Hamiltonian \eqref{e:HNpen}. We restrict the rigorous analysis
to the quadratic interaction potential ($v_{\INTP}(r) = \tfrac{r^2}{2}$),
zero exterior potential ($v_{\EXTP} = 0$), and to the
mass-matrix penalization operator ($\ID-\nusc^2\dl$). The leapfrog scheme is defined as
\begin{equation*}
\left\{\begin{aligned}
p_{n+1/2} &= p_{n} + \frac{\dt}{2} (-\dl) q_{n} \\
q_{n+1} &= q_{n} + \dt \, \penN{M}^{-1} p_{n+1} \\
p_{n+1} &= p_{n+1/2}+ \frac{\dt}{2}  (-\dl)  q_{n+1} \PERIOD
\end{aligned}
\right.
\end{equation*}

Denoting the spectral variables for $k=1,\dots,N-1$
 \begin{equation} \label{e:specvar}
 \left\{\begin{aligned}
\widehat{v}^k & =  \pare{\frac{\delta_k}{1+\nusc^2 \delta_k} }^{1/2} \sqrt{N} P p\\
 \widehat{x}^k & =  \pare{\frac{1+\nusc^2 \delta_k}{\delta_k} }^{1/2} \sqrt{N} P q  \COMMA
 \end{aligned}
 \right.
 \end{equation}
we write
\[
\bmat \widehat{v}_{n+1}^k \\ \widehat{x}_{n+1}^k \emat  = L_k \bmat \widehat{v}_n^k \\ \widehat{x}_n^k \emat \COMMA
\]
where
\[
L_k
=\bmat 1 - \frac{h_k^2}{2} & -h_k + \frac{h_k^3}{4} \\
 h_k & 1 - \frac{h_k^2}{2} \emat
 \COMMA\;\;\;\;\mbox{and}\;\;\;\;h_k = \dt \frac{\delta_k^{1/2}}{\pare{1+\nusc^2 \delta_k}^{1/2}}\PERIOD
\]
Since $\DET\,(L_k)=1$, the standard CFL stability condition is equivalent to \[\abs{\mathrm{Tr}\,(L_k)} \leq 2\]
which is fulfilled  if and only if $h_k \leq 2$ for all $k\leq N-1$.
Thus we arrive at the following bound on the time step
\[
\dt \leq 2 \min_{0\leq k < N} \left(\frac{1+\nusc^2 \delta_k}{\delta_k}\right)^{1/2} \PERIOD
\]

Summarizing the above calculations and recalling \VIZ{eigenval} we have the following characterization
of the stability properties.
\bpro \label{p:stab1}
Suppose $v_{\EXTP} =0$ and consider a harmonic interaction potential $v_{\INTP}(r) = \tfrac{r^2}{2}$
with the mass-matrix penalization $\penN{M}=\ID -\nusc^2 \dl$.
The leapfrog/Verlet integration of the IMMP harmonic Hamiltonian
\eqref{e:HN} is stable in the spectral sense if and only if
\be \label{e:CFL}
\dt \leq \pare{ 4 \nusc^2 + \frac{1}{ \ds N^2 \\ \sin^2\left(\frac{(N-1)\pi}{2N}\right) } }^{1/2}\PERIOD
\ee
\epro

Since we work with a Metropolis correction of the hybrid Monte-Carlo type, we are also interested in the limiting behavior
of the energy variation compared to the temperature, i.e.,
\[
\beta_N (H(p_{n+1},q_{n+1}) - H(p_{n},q_{n})) \COMMA
 \]
 when
 $(p_{n},q_{n})$
are distributed according to the canonical distribution. This quantity gives the average acceptance rate of
the Metropolis correction. The result we present here is similar to \cite{BesStu}
where the authors analyze infinite dimensional sampling
with the standard Metropolis-Hastings Markov chains.
\bpro \label{p:stab2}
Suppose $v_{\EXTP} =0$ and consider a harmonic interaction potential $v_{\INTP}(r) = \tfrac{r^2}{2}$ with
the mass-matrix penalization $\penN{M}=\ID -\nusc^2 \dl$.
Suppose the state variable  $X=(\penN{p},q)$ is a random variable distributed according to the canonical distribution
associated with the mass-matrix penalized Hamiltonian \eqref{e:HN}.
Then the energy variation $\beta_N \Delta H$ after one step
of the leapfrog integration scheme converges in distribution, up to normalization and centering,
to the Gaussian random variable
\[
\frac{\beta_N\Delta H -m_N}{\sigma_N} \xrightarrow[N\to +\infty]{\mathrm{Law}} \mathcal{N}(0,1) \COMMA
\]
with the mean and variance in the infinite size asymptotics for the IMMP method $\nusc > 0 $ and $\dt \equiv\dt_N = o(1)$
\[
  m_N \limop{\sim}_ {N\to +\infty}     \frac{N \dt_N^6}{32 \nusc^{6}} \COMMA \;\;\;\mbox{and}\;\;\;
  \sigma_N^2 \limop{\sim}_ {N\to +\infty} \frac{N\dt_N^6}{16\nusc^{6}} \COMMA
\]
and for the Verlet integration of exact dynamics with $\dt \equiv\dt_N = o(1/N)$
\[
m_N\limop{\sim}_ {N\to +\infty}\frac{5}{8} N^7\dt_N^6  \COMMA\;\;\;\mbox{and}\;\;\;
\sigma_N^2 \limop{\sim}_ {N\to +\infty}  \frac{5}{4}  N^{7} \dt_N^6\PERIOD
\]
\epro
\begin{proof} We start with a canonically distributed state $X=(q,p)$, which is,
by assumption on
the form of the interaction potential,  a Gaussian random vector.
After changing to the spectral coordinates \eqref{e:specvar} we have the spectral representation of the Hamiltonian
\[
\beta_N H = \beta \sum_{k=1}^{N-1}\frac{\delta_k^{1/2}}{2(1+\nusc^2 \delta_k)^{1/2}}\pare{(\widehat{v}^k)^2 + (\widehat{x}^k)^2} \COMMA
\]
and introducing Gaussian random vectors $U$ and $V$ with the identity covariance matrix we can write
\[
\widehat{x}^k = \beta^{-1/2} \frac{(1+\nusc^2 \delta_k)^{1/4}}{\delta_k^{1/4}}U_k\COMMA\;\;\;\mbox{and}\;\;\;
\widehat{v}^k =  \beta^{-1/2} \frac{(1+\nusc^2 \delta_k)^{1/4}}{\delta_k^{1/4}}V_k \PERIOD
\]
We then compute explicitly the change of the Hamiltonian after one step of the leapfrog integration
\begin{equation}\label{deltaH}
\beta_{N} \Delta H = \sum_{k=1}^{N-1} \frac{1}{2} \bmat U_k \\ V_k  \emat^T (L_k^T L_k - \ID ) \bmat U_k \\ V_k  \emat
\PERIOD
\end{equation}
Since $\DET(L_k^T L_k) =1$ the matrix $L_k^T L_k-\ID$ has two positive eigenvalues
$(\lambda_k-1,1/\lambda_k-1)$ which satisfy
\begin{eqnarray*}
&&\lambda_k+1/\lambda_k -2=   \TR(L_k^T L_k-\ID) =\frac{h_k^6}{16} \COMMA \\
&&(\lambda_k-1)^2+(1/\lambda_k -1)^2 =   \TR(L_k^T L_k)^2-2\TR(L_k^T L_k) = \frac{h_k^{12}}{256} + \frac{h_k^6}{8} \PERIOD
\end{eqnarray*}
Combing with \VIZ{deltaH} we find
\[
m_N \equiv \E[\beta_N \Delta H] = \sum_{k=1}^{N-1} \frac{h_k^6}{2^5}\COMMA\;\;\;\mbox{and}\;\;\;
\sigma^2_N \equiv \mathrm{Var}[{\beta_N \Delta H}] = \sum_{k=1}^{N-1}\frac{h_k^6}{2^4}+\frac{h_k^{12}}{2^9}\PERIOD
\]
Moreover, the Lindenberg or simply Lyapunov condition in the general central limit theorem (see \cite{Fel71}) is verified since we work with a sum
of $\chi^2$ random variables, thus concluding the first part of the proof.

Recalling
\[
h_k = \dt \frac{\sin(\frac{k}{N}\frac{\pi}{2})}{(\frac{1}{4N^2}+\nusc^2 \sin^2(\frac{k}{N}\frac{\pi}{2}))^{1/2}} \COMMA
\]
we compute the convergent Riemann sums for $p=6$ and $p=12$.
For the case $\nusc \neq 0$ we have
\[
\lim_{N\to \infty}\frac{1}{N}\sum_{k=1}^{N-1} h^{p}_{k} = \frac{\dt^p}{\nusc^p}\PERIOD
\]
If $\nusc =0$ we obtain
\begin{eqnarray*}
  \lim_{N\to \infty}\frac{1}{N^{p+1}}\sum_{k=1}^{N-1} h^{p}_{k} &=&
  \lim_{N\to \infty}\frac{\dt^p}{N}\sum_{k=1}^{N-1} 2^{p} \sin^{p}\left(\frac{k}{2N}\pi\right) \\
                       &=& \dt^p 2^{p} \int_{0}^{1}\sin^{p}\left(\frac{\pi}{2}x\right)dx\PERIOD
\end{eqnarray*}
Thus for $p=6$ we have that the series sums to $20\dt^6$.
Then the asymptotic behavior follows from the assumption $\dt_N^{12} \ll \dt_{N}^{6}$, and similarly
$\dt_N^{12} N^{12}\ll \dt_{N}^{6} N^6$ in the case $\nusc =0$.
\end{proof}

\brem{\rm
The two propositions proved in this section characterize the restrictions imposed by the
stability of the resulting scheme.
In Proposition~\ref{p:stab1}, stability in the large system size
limit, $N \to + \infty$, is equivalent to the inequality \eqref{e:CFL}. In this case the restriction of the
time-step size is imposed by the numerical integrator. On the other hand
the stability for the scheme which uses a Metropolis corrector is linked to the acceptance
rate of the Metropolis step.
In Proposition~\ref{p:stab1}, stability in the large system size limit is equivalent to
the non-vanishing Metropolis acceptance rate, which is equivalent to bounded from above average energy variation
$m_N$ and bounded variance $\sigma_N$ of  the energy variation.
In either case, the IMMP method ($\nusc>0$) induces a relative increase of order $N$ for the boundary
of numerical stability as compared to the exact dynamics $\nusc=0$ integrated with the Verlet scheme.
}
\erem

\appendix

\section{Surface measures}\label{s:measures}
Let $\mathbb{R}^{d}$ be endowed with the scalar product given by the positive definite matrix $\MASST$, and consider
$\SMANZ{z}$ a family
of sub-manifolds of co-dimension $\nc$ implicitly defined by $\nc$ independent functions
$\SMANZ{z} = \{ x \in \mathcal{M} \SEP \xi_{1}(q)=z_{1},..,\xi_{\nc}(q)=z_{\nc} \}$ for $z$ in a neighborhood of the origin.
For each $z$ in a neighborhood of the origin the \emph{conditional measure} $\delta_{\xi(q)=z}(dq)$ is a measure on $\SMANZ{z}$ defined
in such a way that it satisfies the chain rule for conditional expectations with respect to the Lebesgue measure $dq$,
i.e.,
\be
dq = \delta_{\xi(q) = z}(dq) \, dz\PERIOD
\ee
The {\em surface measure} $\sigma_{T^{*}_{q}\SMANZ{z}}(dp)$  is the Hausdorff measure induced by the metric $\MASST^{-1}$
on the co-tangent space $T^{*}_{q}\SMANZ{z}=\set{ p \SEP \nabla_q^T \xi(q)\MASST^{-1} p  = 0}$; and in the same way,
$\sigma_{\SMANZ{z}}(dq)$  is the  Hausdorff measure induced by the metric
$\MASST$ on the sub-manifold $\SMANZ{z}$. It is important to note that, although this is not explicit in the notation, $\sigma$
is defined with respect to the mass-tensor $\MASST$ of the mechanical system. The Liouville measure $\sigma_{T^{*}\SMANZ{z}}(dp \, dq)$
on the co-tangent bundle
$T^{*}\SMANZ{z}$ is the volume form induced on
\[
T^{*}\SMANZ{z} = \set{(p,q)\SEP \nabla_q^T \xi(q)\MASST^{-1} p  = 0\COMMA\, \xi(q)=z }
\]
by the usual symplectic form $dp \wedge dq$. It can be described in terms of surface measures as follows
\[
\sigma_{T^{*}\SMANZ{z}}(dp \, dq) = \sigma_{T^{*}_q\SMANZ{z}}(dp) \, \sigma_{\SMANZ{z}}(dq) \PERIOD
\]

Finally, the co-area formula (see \cite{Fed69} for a general reference) defines the relative probability density between
$\delta_{\xi(q)=z}(dq)$ and $\sigma_{\SMANZ{z}}(dq)$.
\begin{Pro}[Co-area formula]
Given the invertible Gram matrix associated with the constraints $\xi(q) =z$ in a neighborhood of $\SMANZ{z} = \{ q \SEP \xi(q)= z \}$
\[
G(q) = \nabla_q^T \xi\, \MASST^{-1} \nabla_q \xi\COMMA
\]
one has
\[
\delta_{\xi(q)=z}(dq) =  \frac{1}{\sqrt{ \DET G(q)}} \sigma_{\SMANZ{z}}(dq) \PERIOD
\]
\end{Pro}
\section{Langevin processes}\label{s:langevinapp}
Defining the Poisson bracket
\[
\poisson{\ph_1,\ph_2} = \Dp^T\ph_1  \Dq \ph_2  - \Dp^T\ph_2  \Dq \ph_1 \COMMA
\]
and the dissipation tensor
\[
\diss(q) = \sigma q \COMMA
\]
where $\sigma$ is the fluctuation matrix in Definition~\ref{d:langevin},
the Markov generator of the Langevin process in Definition~\ref{d:langevin} is
\[
\mathcal{L} = \poisson{\,\cdot\, , H} + \frac{1}{\beta}
\poisson{\diss,\poisson{\diss^T, \, \cdot \, }\EXP{- \beta H}}\EXP{\beta H} \PERIOD
\]
The generator $\mathcal{L}$ satisfies
\[
\int \ph_1 \,\mathcal{L} (\ph_2) \EXP{- \beta H} \,dp \, dq =
\int \mathcal{L}^* (\ph_1) \,\ph_2 \,\EXP{- \beta H} \,dp \, dq \COMMA
\]
where
\[
\mathcal{L}^* = \poisson{\,\cdot\, , -H} + \frac{1}{\beta}
\poisson{\diss,\poisson{\diss^T, \, \cdot \, }\EXP{- \beta H}}\EXP{\beta H} \PERIOD
\]
The generator $\mathcal{L}^*$ defines a Langevin process with the time-reversed Hamiltonian ($-H$).
Reversibility of the process implies that the canonical measure is stationary.
Furthermore, if the initial state of the system is a canonically distributed random variable, the probability distribution
of a trajectory after the time-reversal is given by a Langevin process with the generator $\mathcal{L}^*$.
When $H$ has the form $H(p,q) =  \tfrac{1}{2} p^{T} \MASST^{-1} p    + V(q) $, reversal of impulses ($p\to -p$) leads
to time-reversed dynamics, and a process with
generator $\mathcal{L}^*$ can be constructed by the following simple steps:
 \begin{enumerate}
 \item Reverse momenta ($p\to -p$).
 \item Draw a random path with generator $\mathcal{L}$.
 \item Reverse again momenta ($p\to -p$).
 \end{enumerate}
When holonomic constraints, for instance, of the form
\[
\Xi(p,q) = \zeta \Leftrightarrow \syst{
p^{T} M^{-1} \nabla_{q} \xi &=& 0 \\
\xi(q) &=& z
}
\]
are introduced, it is useful to define the Poisson bracket on the co-tangent bundle $T^*\SMANZ{z}$
\[
\poisson{\ph_1,\ph_2}_{\SMANZ{z}} = \poisson{\ph_1,\ph_2} - {\ds \sum_{a,b}}\poisson{\ph_1,\Xi^a} \Gamma^{-1}_{a,b} \poisson{\Xi^b,\ph_2} \COMMA
\]
where $\Gamma$ is the symplectic Gram matrix of the full constraints
\[
\Gamma^{a,b} = \poisson{\Xi^a,\Xi^b} \PERIOD
\]
As a basic result of symplectic geometry (see \cite{Arn89}), one recovers the divergence formula with respect to the bracket
$\poisson{\,\cdot \,\,, \cdot\,}_{\SMANZ{z}} $ and the Liouville measure $\sigma_{T^{*}\SMANZ{z}}(dp \, dq)$
\[
\int \poisson{\,\cdot \,, \cdot\,}_{\SMANZ{z}} \sigma_{T^{*}\SMANZ{z}}(dp \, dq) = 0 \PERIOD
\]
Given a constrained Langevin process in a stochastic differential equation form
\begin{eqnarray*}
    \dot{q} &=& \Dp H \COMMA \\
    \dot{p} &=& -\Dq H -\gamma \dot{q} + \sigma \dot{W} - \Dq \xi \,\dot{\lambda}\COMMA
   \end{eqnarray*}
where $\lambda$ are Lagrange multipliers associated with the constraints $\xi(q)=0$,
adapted with respect to the noise $\dot{W}$, the process $\{p_t,q_t\}_{t\geq 0}$ obeys hidden
velocity constraints and is characterized
by the stochastic differential equations
\begin{eqnarray*}
    \dot{q} &=& \Dp H  +  \Dp \Xi \,\dot{\Lambda}\COMMA \\
    \dot{p} &=& -\Dq H -\gamma \dot{q} + \sigma \dot{W} - \Dq \Xi\,\dot{\Lambda}\COMMA
\end{eqnarray*}
where $\Lambda$ are Lagrange multipliers associated with the full constraints $\Xi(p,q)=0$.
The Markov generator of this process can be written in the form
\[
\mathcal{L}_{\SMANZ{z}} = \poisson{\,\cdot \,, H}_{\SMANZ{z}} + \frac{1}{\beta}
\poisson{\diss,\poisson{\diss^T, \, \cdot \, }_{\SMANZ{z}}
\EXP{- \beta H}}_{\SMANZ{z}}\EXP{\beta H} \COMMA
\]
demonstrating the reversibility with respect to the constrained canonical measure
$\EXP{-\beta H} \sigma_{T^{*}\SMANZ{z}}(dp \, dq)$.

\section{Exact sampling of fluctuation/dissipation perturbations}\label{s:exactflucdiss}
In this section, we recall how to perform exact sampling of fluctuation/dissipation perturbations.
Since we only work with impulses, we refer to the system by using the impulse variables $p$ only.
Note that throughout the paper, we also use extended variables $(p,p_z)$,
however, the presentation that follows covers general cases.
The kinetic energy of the system is $\tfrac{1}{2}  p^{T}\MASST  p$. We impose constraints
$p^{T} \MASST^{-1} \Dq \xi = 0$ on impulses, thus $p\in T^{*}_{q}\mathcal{M}$
and hence the associated orthogonal projector on $T^{*}_{q}\mathcal{M}$ is
\[
P  = \ID - \Dq \xi \,G^{-1} \Dq^T \xi \,\MASST^{-1} \PERIOD
\]
The stochastic differential equations of motion on impulses that are integrated on a time-step interval are
\begin{equation}\label{e:sdeimp}\syst{
 &\dot{p}= -\gamma \MASST^{-1}p + \sigma \dot{W} - \Dq \xi \, \dot{\lambda}\COMMA &\\
 &p^{T} \MASST^{-1} \Dq \xi = 0\COMMA & \TAG{C_p}
}\end{equation}
with the usual fluctuation/dissipation relation $\sigma \sigma^T = 2\beta^{-1}\gamma$.
The Gaussian distribution of impulses
\be\label{e:maxwell}
\frac{1}{Z} \EXP{-\frac{\beta}{2} p^{T}\MASST^{-1}  p} \sigma_{T^{*}_{q}\mathcal{M}}(dp)
\ee
is invariant under the dynamics \VIZ{e:sdeimp}.
\bpro[Exact sampling of stochastic perturbation]\label{p:exactscheme}
Given the mass matrix $M$, suppose either $\dt$ or $\gamma$  are small enough so that the condition
\be
\frac{\delta t}{2}  \MASST^{-1}   \leq \gamma
\ee
holds in the sense of symmetric semi-definite matrices. Let $U$ be a centered and normalized Gaussian vector.
Consider the mid-point Euler scheme with constraints
\begin{equation} \label{e:midpoint}\syst{
& p_{n+1}= p_{n}  - \frac{\dt}{2} \gamma \MASST^{-1}(p_{n} + p_{n+1)} + \sqrt{\dt} \sigma U - \Dq \xi \,\lambda_{n+1} &\\
& p_{n+1}^T\MASST^{-1}  \Dq \xi = 0 \COMMA & \TAG{C_p} %
}\end{equation}
where $\lambda_{n+1}$ is the Lagrange multiplier associated with the constraint $(C_p)$.
The Markov kernel defined by the transition $p_{n}\to p_{n+1}$ is reversible with respect to
the Gaussian distribution \eqref{e:maxwell}.
\epro
\begin{proof}
After calculating the Lagrange multiplier the expression \VIZ{e:midpoint} can be written as
\[
p_{n+1} = p_{n}  - \frac{\dt}{2} P \gamma P^{T} \MASST^{-1}(p_{n} + p_{n+1}) + \sqrt{\dt} P \sigma U \PERIOD
\]
Consider the new variable
$\tilde{p} = \beta^{1/2} \MASST^{-1/2} p$,
and define the symmetric matrix
$$
\ZETAM \equiv \frac{\dt}{2} \MASST^{-1/2} P \gamma P^{T} \MASST^{-1/2} \COMMA
$$
as well as $\KAPPAM$, such that $\KAPPAM \KAPPAM^{T} = \ZETAM$.
In terms of these new variables we obtain from (\ref{e:midpoint})
\be \label{e:midpoint_tilde}
\tilde{p}_{n+1} = (\ID + \ZETAM)^{-1}(\ID-\ZETAM)\,\tilde{p}_n + 2 (\ID+\ZETAM)^{-1}\KAPPAM \,U \PERIOD
\ee
Moreover, the product measure $\sigma_{T^*_q\SMAN}(dp_n)\,\sigma_{T^*_q\SMAN}(dp_{n+1})$ is the measure induced on the linear subspace
of constraints by the scalar product $\MASST^{-1}$ and the Lebesgue measure $dp_n \, dp_{n+1}$.
Thus in the variables $( \tilde{p}_n, \tilde{p}_{n+1})$ this measure becomes, up to a constant, the measure induced by the usual Euclidean structure.
As a consequence the $\log$ density of the random variable $(\tilde{p}_n, \tilde{p}_{n+1})$ defined by
\eqref{e:midpoint_tilde} with respect to this latter measure is equal to
\begin{eqnarray*}
 && - \frac{1}{2} \abs{\tilde{p}_n}^2 - \frac{1}{8}\left(\tilde{p}_{n+1}- (\ID+\ZETAM)^{-1}(\ID-\ZETAM)\,\tilde{p}_n\right)^T
                  \ZETAM^{-1}(\ID+\ZETAM)^2 \left(\tilde{p}_{n+1}- (\ID+\ZETAM)^{-1}(\ID-\ZETAM)\,\tilde{p}_n \right) \\
&=& - \frac{1}{8}\tilde{p}_{n+1}^T \ZETAM^{-1}(\ID+\ZETAM)^2\,\tilde{p}_{n+1}
    - \frac{1}{8}\tilde{p}_{n}^T \ZETAM^{-1}(\ID+\ZETAM)^2\,\tilde{p}_{n} \COMMA
\end{eqnarray*}
which is indeed symmetric between $\tilde{p}_{n}$ and $\tilde{p}_{n+1}$. Hence
we have shown the reversibility of the induced Markov kernel and consequently
stationarity of the canonical Gaussian distribution.
\end{proof}

\bibliographystyle{plain}

\end{document}